\documentclass[final,3p,10pt]{elsarticle}

\usepackage{enumitem}
  \setlist{nosep} %

\usepackage[makeroom]{cancel}
\usepackage{amssymb,amsthm,mathtools,xfrac}
\usepackage{graphicx}
\usepackage{stmaryrd}
\usepackage{verbatim}
\usepackage{soul}
\usepackage{comment}
\usepackage{color}
\usepackage{cleveref}
\usepackage{subcaption}
\renewcommand{\d}{\textnormal{d}}

\usepackage{xcolor}

\newcommand{\tcb}[1]{\textcolor{black}{#1}}

\usepackage{siunitx}
\DeclareSIUnit\str{str} %
\DeclareSIUnit\sh{sh} %
\DeclareSIUnit\erg{erg} %

\usepackage{algorithm,algorithmicx}
\usepackage[noend]{algpseudocode}

\DeclareMathOperator*{\argmin}{arg\,min}

\journal{Journal of Computational Physics}

\begin{document}

\begin{frontmatter}

\title{Moment-based adaptive time integration for thermal radiation transport}
\author[LANL]{Ben S. Southworth}
\cortext[mycorrespondingauthor]{Corresponding author}
\ead{southworth@lanl.gov}
\author[LANL]{Steven Walton}
\author[LLNL]{Steven B. Roberts}
\author[LANL]{HyeongKae Park}

\address[LANL]{Theoretical Division, Los Alamos National Laboratory, P.O. Box 1663, Los Alamos, NM 87545 U.S.}
\address[LLNL]{Center for Applied Scientific Computing, Lawrence Livermore National Laboratory, 7000 East Ave, Livermore, CA 94550 U.S.}

\begin{abstract}
In this paper we develop a framework for moment-based adaptive time integration of deterministic multifrequency thermal radiation transpot (TRT). We generalize our recent semi-implicit-explicit (IMEX) integration framework for gray TRT to multifrequency TRT, and also introduce a semi-implicit variation that facilitates higher-order integration of TRT, where each stage is implicit in all components except opacities. To appeal to the broad literature on adaptivity with Runge--Kutta methods, we derive new embedded methods for four asymptotic preserving IMEX Runge--Kutta schemes we have found to be robust in our previous work on TRT and radiation hydrodynamics. We then use a moment-based high-order-low-order representation of the transport equations. Due to the high dimensionality, memory is always a concern in simulating TRT. We form error estimates and adaptivity in time purely based on temperature and radiation energy, for a trivial overhead in computational cost and memory usage compared with the base second order integrators. We then test the adaptivity in time on the tophat and Larsen problem, demonstrating the ability of the adaptive algorithm to naturally vary the timestep across 4--5 orders of magnitude, ranging from the dynamical timescales of the streaming regime to the thick diffusion limit. 
\end{abstract}
\begin{highlights}
\item Adaptive implicit-explicit integration developed for thermal radiation transport
\item Computationally cheap and simple to apply
\item Naturally follows dynamical timescale across 4--5 orders of magnitude
\end{highlights}


\end{frontmatter}

\section{Introduction}

Consider the multifrequency time-dependent thermal radiation transport (TRT) and material temperature equations without scattering given by 
\begin{subequations}\label{eq:ho}
\begin{align}\label{eq:ho-I}
  \frac{1}{c} \frac{\partial I}{\partial t} & = - \Omega \cdot \nabla I - \sigma(\nu,T) I + \sigma(\nu,T)B(\nu,T),\\
  \rho C_v \frac{\partial T}{\partial t} & = \int_0^\infty \int_{S^d} \Big( \sigma(\nu,T) I  - \sigma(\nu,T)B(\nu,T) \Big)\d\Omega\d\nu .\label{eq:ho-T}
\end{align}
\end{subequations}
Here, $I(\mathbf{x},\Omega,\nu,t)$ is the specific radiation intensity for photons traveling in direction $\Omega$ with frequency $\nu$, $T(\mathbf{x},t)$ is the material temperature, $c$ is the speed of light, $C_v$ is the specific heat (which in certain cases may depend on $T$), $\rho$ is the material density, and $\sigma(\nu,T)$ are energy- and temperature-dependent material opacities. In complex physics regimes, opacities do not have a known analytical form, and are defined via a tabular data look up based on temperature values. The emission source $B(\nu,T)$ is given as the Planckian function
\begin{equation*}
    B(\nu,T) = \frac{2h\nu^3}{c^2}\left( e^{h\nu/kT} - 1\right)^{-1},
\end{equation*}
where $h$ and $K$ are the Planck and Boltzmann constants, respectively \cite{Mihalas.1986}.

In non-relativistic regimes, we typically want to step well over the transport and collisional timescales, and as a result \eqref{eq:ho} is treated implicitly in time. Due to the high dimensionality and stiff coupling of variables and dimensions, though, implicit integration of \eqref{eq:ho} is very expensive. There are also multiple forms of stiffness, some of which is physical and must be respected, e.g. to capture the correct wavefront propagation, and other that is numerical, wherein we can step well over the numerical CFL without loss of physics fidelity. Additionally, the evolution of photons and coupling to material temperature can experience a wide range of dynamic timescales in a single simulation. As radiation initially enters an optically thin domain, radiation energy may evolve on the transport/advective time scale, while the evolution of the radiation energy density in optically thick regions happens on a diffusion timescale, which is a function of material opacities and specific heat. Changes in material or corners around which energy and temperature must propagate can also lead to local physical stiffness. Altogether, the dynamical timescale arising in the simulation of radiation flow in heterogeneous media evolves with the problem, and can span a wide range.

Altogether, these properties suggest adaptive time stepping to be a natural topic of interest in TRT simulation. Effectively, we want to take the maximum possible timestep (to minimize the number of high-dimensional implicit solves), while also ensuring we maintain physics fidelity.
For general systems of ordinary differential equations (ODEs), there is an extensive body of literature on adaptively selecting a timestep based on a local temporal error estimate, e.g. \cite{Ceschino.1961,Gustafsson.1994,Soderlind.2002}. Local error estimates for one-step time integration methods typically come from comparing a primary solution to an embedded solution of a different order. While it is relatively common for Runge--Kutta methods from the literature to be equipped with embeddings, it is significantly less common for more sophisticated integrators including the class of partitioned Runge--Kutta (PRK) methods considered in this work. In fact, our previous work has found it necessary to use the class of asymptotic preserving IMEX schemes for hyperbolic equations with stiff relaxation terms as defined in \cite{Pareschi.2005} for accurate and stable IMEX integration applied to TRT \cite{imex-trt} and radiation hydrodynamics \cite{rad-hydro}. To the best of our knowledge, we are aware of no such schemes with derived embedded methods in the literature. 

Interestingly, there has been very little work to date in the literature on adaptivity in time for TRT, and the large codes developed at U.S. Department of Energy labs tend to use ad-hoc measures based on relative change in solution, or simply follow the sound speed CFL determined by a coupled, e.g., hydrodynamics, code. We hypothesize one reason for this is that typically adaptive methods are built based on comparing the new solution constructed from an order-$p$ and order-$p+1$ integrator, and the majority of deterministic transport codes use a first -order semi-implicit Euler method, where
\begin{equation}\label{eq:limex-o1}
    y_{n+1} = y_n + \Delta t \mathcal{N}(y_n, y_{n+1}),
\end{equation}
with opacities evaluated at the previous time step \cite{Larsen.1988}. Opacities $\sigma$ typically depend nonlinearly on frequency and temperature, and in many realistic problems are represented via tabular data profiles, each of these aspects making it difficult to resolve opacity dependence implicitly. Due to the nonlinear coupling of opacities with other components of the system, this makes the use of general higher-order integrators challenging. Higher-order DIRK methods are considered in \cite{maginot2016high}, but here opacities must either be constant or resolved implicitly, which are not always realistic or practical requirements, respectively. Predictor-corrector schemes have been considered in, e.g. \cite{Park.2013}, but without adaptivity, and such schemes typically have relatively poor stability as observed in the order reduction for multifrequency problems seen in \cite{Park.2013}. Mixed-order schemes are also considered in \cite{ghassemi2020multilevel}, where low-order moments are integrated with higher order schemes to improve accuracy without significantly increasing computational cost, but adaptivity is not considered. 

In this paper we develop a framework for adaptive time integration of deterministic multifrequency TRT. We generalize our recent semi-implicit-explicit integration framework for gray TRT \cite{imex-trt} to multifrequency, and also introduce a semi-implicit variation (higher-order integration, where each stage is implicit in all components except opacities). This allows us to use PRK methods and the broad literature on adaptivity with Runge--Kutta methods, in the context of the semi-implicit framework of \cite{Boscarino.2015,Boscarino.2016} that facilitates nonlinear partitioning. To facilitate this, we derive embedded methods for four asymptotic preserving IMEX Runge--Kutta schemes we have found to be robust in our previous work on TRT \cite{imex-trt} and radiation hydrodynamics \cite{rad-hydro}. We then use a moment-based representation of the transport equations \cite{gol1964quasi,Anistratov.1993}, where we integrate the ``low-order'' (LO) zeroth and first moment of \eqref{eq:ho} coupled to the ``high-order'' (HO) transport equation of specific intensity \eqref{eq:ho}. Often the LO equations are used purely as a tool to accelerate the implicit solution of the HO equations. Here we treat the LO moments as auxiliary integrated quantities in time in order to (i) use them to form a partitioning as in \cite{imex-trt}, and (ii) as the basis for our error estimates and adaptivity. Due to the high dimensionality, memory is always a concern in simulating TRT. We form error estimates and adaptivity in time purely based on temperature and radiation energy, for a trivial overhead computational cost and memory usage compared with the base second order integrator. Specifically, using temperature and radiation energy allow us to incorporate adaptivity without additional storage and operations on the full specific intensity space. 

In \Cref{sec:holo} we present the semi-implicit and semi-implicit-explicit nonlinear partitioning of TRT that we use in the context of semi-implicit Runge--Kutta methods, which facilitate using arbitrary order partitioned Runge--Kutta methods and their corresponding error estimators. In \Cref{sec:time} we review semi-implicit Runge--Kutta methods and adaptive Runge--Kutta methods. We then choose four specific second-order partitioned integrators we have found to be robust for TRT problems, and derive new embedded error estimators with certain optimal properties for these methods, and present the framework for moment-based adaptive time stepping. Note, the framework developed here applies to general moment-based approaches to TRT. We specifically use high-order-low-order (HOLO) moment discretizations \cite{Park.2012,Park.2013,Park.2020} because in our testing they have demonstrated better IMEX stability for transport, and they also provide a natural/consistent collapse to a gray moment system that we use for adaptivity and IMEX partitioning. The new framework is then applied to the gray tophat and multifrequency Larsen problem in \Cref{sec:results}, demonstrating the ability to naturally change the timestep to respect the dynamic timescale across many orders of magnitude in a single simulation.

\section{TRT moment equations and semi-implicit(-explicit) formulation}\label{sec:holo}

Here we consider deterministic S$_N$ transport discretizations, where the angular integral operator is discretized using discrete ordinates \cite{Adams.2002} and the frequency integral approximated using a multigroup energy discretization. We use a linear upwind discontinuous Galerkin discretization in space, with mass lumping to help maintain positivity \cite{Park.2020}. We begin by introducing a general Runge--Kutta framework for (semi-)implicit-explicit integration in \Cref{sec:holo:simex}, and proceed to introduce HOLO-based implicit-explicit and semi-implicit integration methodologies for TRT in \Cref{sec:imex:trans1}, following the gray formulation from \cite{imex-trt}. 

\subsection{Semi-implicit(-explicit) Runge Kutta}\label{sec:holo:simex}

Post spatial discretization, TRT results in a nonlinear autonomous set of ODEs in time, with general form
\begin{align}\label{eq:ode}
    \frac{\partial y}{\partial t} = \mathcal{N}(y),
    \qquad y \colon \mathbb{R} \to \mathbb{R}^{N}
\end{align}
Runge-Kutta (RK) methods \cite{butcher1996history} approximate the integration of \eqref{eq:ode} in time using a combination of Euler-like stage steps/solves. We adopt the standard notation of a Butcher tableaux,
$
\begin{array}
{c|c}
\mathbf{c} & A\\
\hline
&\mathbf{b}^T
\end{array},
$
where $A$ corresponds to RK coefficients, $\mathbf{b}$ to weights, and $\mathbf{c}$ to quadrature points within a time-step. A generic $s$-stage RK method takes the form
\begin{align*}
    Y_{(i)} &= y^n + \Delta t \sum_{j=1}^s a_{ij} \mathcal{N}(Y_{(j)})\hspace{4ex}\text{for }i=1,\dots,s,\\
    y^{n+1} &= y^n + \Delta t \sum_{j=1}^s b_j \mathcal{N}(Y_{(j)}).
\end{align*}
Strictly lower triangular $A$ yields an explicit method, and lower triangular yields a diagonally implicit RK (DIRK) method (we do not consider fully implicit methods here, which yield more complex implicit stage equations). 

As discussed previously, for realistic TRT simulations we typically do not want to consider the problem as purely implicit due to the nonlinear dependence on material opacities. Significant work has been done on partitioned and additive RK (ARK) methods \cite{Ascher.1997,Kennedy.2003tv4}, which are designed to integrate equations with a natural additive partition into stiff and nonstiff parts, $\mathcal{N}(y) = \mathcal{N}_E(y) + \mathcal{N}_I(y)$, where we treat $\mathcal{N}_E$ explicitly and $\mathcal{N}_I$ implicitly. However, in TRT the opacities that we want to treat explicitly (or effectively linearize) are nonlinearly coupled to the stiff components that we must treat implicitly.

To address the nonlinear coupling arising in opacities, we follow the semi-implicit strategy of \cite{Boscarino.2015,Boscarino.2016}: we introduce an auxiliary $*$-variable and apply IMEX-RK methods to solve the resulting coupled equations:
\begin{equation}\label{eq:mod-ode}
    \frac{\partial y^*}{\partial t} = \mathcal{N}(y^*,y), \hspace{3ex}
    \frac{\partial y}{\partial t} = \mathcal{N}(y^*,y).
\end{equation}
Note that because the right-hand sides are identical, $\mathcal{N}(y^*,y)$, exact integration will yield identical solutions, $y(t) = y^*(t)$. By duplicating the equations, however, we can partition our equations such that stiff terms are evaluated at $y$ and non-stiff terms at $y^*$; we then proceed to apply partitioned integrators that are explicit in $y^*$ and implicit in $y$ to the coupled system \eqref{eq:mod-ode}, regardless of nonlinear coupling between the two variables. This will allow us to nonlinearly partition opacities, e.g., $\sigma(\nu,T^*)I$, where this term will be treated explicitly in $T$ and implicitly in $I$. Thus define an IMEX-RK Butcher tableaux
\begin{equation}\label{eq:imex-tableaux}
\renewcommand\arraystretch{1.2}
\begin{array}
{c|c}
\tilde{\mathbf{c}} & \tilde{A}\\
\hline
& \tilde{\mathbf{b}}^T
\end{array}
,\hspace{3ex}
\begin{array}
{c|c}
\mathbf{c} & A\\
\hline
&\mathbf{b}^T
\end{array}
\end{equation}
where tilde-coefficients denote an explicit scheme and non-tilde coefficients represent a DIRK scheme. We assume that $\mathbf{b} = \tilde{\mathbf{b}}$, which simplifies the resulting semi-implicit(-explicit) RK scheme applied to take the form
\begin{subequations}\label{eq:limex}
\begin{align}\label{eq:limex-exp}
    Y^*_{(i)} &= y_n + \Delta t \sum_{j=1}^{i-1} \tilde{a}_{ij} \mathcal{N}(Y^*_{(j)},Y_{(j)}),\\
    Y_{(i)}   &= y_n + \Delta t \sum_{j=1}^i {a_{ij}} \mathcal{N}(Y^*_{(j)},Y_{(j)}),\label{eq:limex-imp} \\
    y_{n+1} &= y_n + \Delta t \sum_{j=1}^s {b_{j}} \mathcal{N}(Y^*_{(j)},Y_{(j)}).\label{eq:limex-sol}
\end{align}
\end{subequations}

\subsection{Continuous TRT moment equations}

Moment-based methods typically take the first two moments of \eqref{eq:ho} in momentum (equivalently the first two moments in angle, along with integrating over energy), yielding the total radiation energy and radiation flux
\begin{equation}\label{eq:moments}
    E \coloneqq \frac{1}{c}\int \int I \d\nu\d\boldsymbol{\Omega}, \hspace{5ex} 
    \mathbf{F} \coloneqq \int \int I\boldsymbol{\Omega} \d\nu\d\boldsymbol{\Omega},
\end{equation}
and then introduce some closure to the radiation flux equation that arises when integrating the right-hand side of \eqref{eq:ho-I}. Here we use the HOLO class of methods \cite{Park.2012,Park.2013,Park.2020} which introduce a $\gamma$-consistency term to the radiation flux equation, wherein the resulting LO system takes the form
\begin{subequations}\label{eq:lo}
\begin{align}
    \frac{\partial E}{\partial t} &= -\nabla \cdot \mathbf{F} - c\sigma_E(T)E + ac\sigma_P(T) T^4,\\
    \frac{1}{c}\frac{\partial \mathbf{F}}{\partial t} & = -\frac{c}{3}\nabla E - \sigma_R(T) \mathbf F + \gamma(I,T)cE, \label{eq:moment-eqs-F}\\
    \rho c_v \frac{\partial T}{\partial t} & = -ac\sigma_P(T)T^4 + c\sigma_E(T) E.
\end{align}
\end{subequations}
Here, $\sigma_E,\sigma_P$, and $\sigma_R$ are the radiation-energy-, Planck-, and Rosseland-weighted opacities, respectively, given by 
\begin{equation}\label{eq:collapse}
    \sigma_E(I,T) \coloneqq \frac{\int\int \sigma I \d\nu\d\boldsymbol{\Omega}}{\int \int I \d\nu\d\boldsymbol{\Omega}}, \hspace{4ex}
    \sigma_P(T) \coloneqq \frac{\int \sigma B \d\nu}{\int B \d\nu}, \hspace{4ex}
    \sigma_R(T) \coloneqq \frac{\int \frac{\partial B}{\partial T} \d\nu}{\int \frac{1}{\sigma}\frac{\partial B}{\partial T} \d\nu},
\end{equation}

Note that the moment equations in \eqref{eq:lo} are auxiliary in the sense that the transport equation in \eqref{eq:ho} contains all of the information of the moment equations. A typical approach of moment-based methods is to apply the first-order semi-implicit discretization in time \eqref{eq:limex-o1}, where opacities are linearized and all other terms treated implicitly, and proceed to use the moment equations as a nonlinear preconditioner for the high-order implicit system. First, a ``transport sweep'' is applied, which solves the discrete implicit form of \eqref{eq:ho-I} for all discretized angles and energies, with fixed temperature and opacities. Then, the collapsed opacities \eqref{eq:collapse} and moment closure (in this case $\gamma$ \eqref{eq:lo}) are formed from the updated angular intensity, and the nonlinear LO equations are solved for updated $\{E,\mathbf{F}, T\}$ with fixed opacities and closure. This process is repeated until convergence of the implicit system.

A benefit of moment-based methods is the ability to choose different spatial discretizations for the HO and LO systems, each of which may facilitate different desirable properties in terms of discretization accuracy, linear/nonlinear solver efficiency, or ease of multiphysics coupling. Finite volume (FV) methods are a natural choice for the LO discretization, as they are widely used in lab codes and also make for natural coupling to FV-based hydrodynamics codes, but care must be taken such that the LO FV discretization satisfies the thick diffusion limit. To that end, we use the HOLO subcell FV discretization of the LO equations as introduced in \cite{Park.2020}. The time integration framework is also largely agnostic to discretization of the integrals in \eqref{eq:moments} and \eqref{eq:collapse}; here we use discrete ordinates in angles and a standard multigroup approximation in energy. 

\subsection{Semi-implicit(-explicit) moment formulation}\label{sec:imex:trans1}

For embedded adaptivity, one needs a base solution method with order $p>1$, so that an embedded solution of order $p-1\geq 1$ can also be formed. This makes the classical first-order semi-implicit method \eqref{eq:limex-o1} widely used in TRT not a viable option. Here we generalize the first-order semi-implicit method to higher orders via the framework of \cite{Boscarino.2015,Boscarino.2016} (see \Cref{sec:holo:simex}). Note that memory is critical in TRT simulations due to the high dimensionality of the problem, and copying all state variables to have both explicit and implicit stage vectors at first glance seems problematic. We want to emphasize that in both of the formulations that follow and the integration schemes introduced in \Cref{sec:time}, the only state variable for which an explicit copy must be stored is temperature, which does not have angular and energy dependence, so the auxiliary variable formulation is a marginal additional memory cost. \tcb{Detailed algorithms and pseudocode for semi-implicit-explicit integration in the context of TRT can be found in \cite[Sec. 3.2 and 4]{imex-trt}.}

\subsubsection{Semi-implicit TRT}

Broadly we want to treat the stiff differential operators and asymptotic behavior of TRT implicitly, but do not want to couple opacity data nonlinearly into our implicit solve. As described previously, this is typically accomplished by linearizing opacities, density, and heat capacity at the previous solution, and treating all other terms in \eqref{eq:ho} implicitly. In the case of HOLO methods, we augment \eqref{eq:ho} with the low-order moment equations \eqref{eq:lo} and alternate between solving the two systems until they are consistent up to a desired tolerance (that is, the discrete moments $E$ and $\mathbf{F}$ match the integrated moments \eqref{eq:moments} from $I$). Often this is formalized post time discretization \cite{Park.2013,Park.2020}, but here we consider a HOLO formulation pre time discretization to allow using arbitrary partitioned semi implicit integrators. To do this, we first need to define the concept of \emph{partitioned collapsed opacities}, which use two distinct temperatures $\{T*,T\}$ for opacities and weight-functions: 
\begin{equation*}%
    \sigma_E(T^*,I) \coloneqq \frac{\int\int \sigma(\nu,T^*) I \d\nu\d\boldsymbol{\Omega}}{\int \int I \d\nu\d\boldsymbol{\Omega}}, \hspace{2ex}
    \sigma_P(T^*,T) \coloneqq \frac{\int \sigma(\nu,T^*) B(\nu,T) \d\nu}{\int B(\nu,T) \d\nu}, \hspace{2ex}
    \sigma_R(T^*,T) \coloneqq \frac{\int \frac{\partial B(\nu,T)}{\partial T} \d\nu}{\int \frac{1}{\sigma(\nu,T^*)}\frac{\partial B(\nu,T)}{\partial T} \d\nu},
\end{equation*}
where in the spirit of \eqref{eq:mod-ode} we will treat $*$-variables explicitly in our integration scheme and non-$^*$-variables implicitly. We then arrive at a complete nonlinearly partitioned set of equations given by
\begin{subequations}\label{eq:semi-trt}
\begin{align}
    \frac{\partial I}{\partial t}  &= - c\boldsymbol{\Omega} \cdot \nabla I - c\sigma(\nu,T^*) I + c\sigma(\nu,T^*)B(\nu,T) &&\coloneqq \mathbf{N}_I^{\textnormal{semi}}(T^*,I,T),\label{eq:semi-I}\\
\frac{\partial E}{\partial t} &=-\nabla \cdot \mathbf{F} - c\sigma_E(T^*,I)E + ac\sigma_P(T^*,T)T^4&&\coloneqq \mathbf{N}_E^{\textnormal{semi}}(T^*,I,E,\mathbf{F},T),\\
    \frac{\partial \mathbf{F}}{\partial t} &=- \frac{c^2}{3}\nabla E - c\sigma_R(T^*,T)\mathbf F   + \gamma(I,T^*) &&\coloneqq \mathbf{N}_{\mathbf{F}}^{\textnormal{semi}}(T^*,I,E,\mathbf{F},T),\\
    \frac{\partial T}{\partial t} &= - \frac{1}{\rho c_v}\left(ac\sigma_P(T^*,T)T^4 - c\sigma_E(T^*,I)E\right) &&\coloneqq \mathbf{N}_T^{\textnormal{semi}}(T^*,I,E,T).
\end{align}
\end{subequations}
In practice, we solve each semi-implicit system that arises in integrating \eqref{eq:semi-trt} via standard HOLO iterations, first solving $I$ for fixed $T$, and then solving the LO equations over $\{E,\mathbf{F},T\}$ given current specific intensity $I$, and repeating until convergence of the HO and LO systems. Note, in the inner LO solve, the $T$ dependence in opacities and $\gamma$ is fixed in the inner nonlinear solve; i.e., the only nonlinearity in $T$ that is resolved in the inner LO solve is the $T^4$ term. More broadly, HO opacities $\sigma(\nu,T^*)$ in \eqref{eq:semi-I} are linearized about $T$ from the previous time step or stage, but collapsed opacities in the LO equations are updated each outer nonlinear HOLO iteration based on updated $I$ and $T$ values.

\subsubsection{Semi-implicit-explicit TRT}

A single implicit transport solve is computationally very expensive due to the high dimensionality, ill-conditioning, nonlinearities, and advective nature of the equations. Recently, we proposed a new one-sweep implicit-explicit time integration scheme for gray TRT based on capturing the asymptotic streaming and thick diffusion regimes implicitly, and treating the remaining equations explicitly. Here we extend the IMEX approach to the multifrequency setting in an analogous way by breaking part of the coupling between the HO and LO equations. We do this by treating the coupling of $B(\nu,T)$ in \eqref{eq:semi-I} expicitly in $T$:
\begin{subequations}\label{eq:imex-trt}
\begin{align}
    \frac{\partial I}{\partial t}  &= - c\boldsymbol{\Omega} \cdot \nabla I - c\sigma(\nu,T^*) I + c\sigma(\nu,T^*)B(\nu,T^*) &&\coloneqq \mathbf{N}_I^{\textnormal{imex}}(T^*,I),\\
\frac{\partial E}{\partial t} &=-\nabla \cdot \mathbf{F} - c\sigma_E(T^*,I)E + ac\sigma_P(T^*)T^4&&\coloneqq \mathbf{N}_E^{\textnormal{imex}}(T^*,I,E,\mathbf{F},T),\\
    \frac{\partial \mathbf{F}}{\partial t} &=- \frac{c^2}{3}\nabla E - c\sigma_R(T^*)\mathbf F   + \gamma(I,T^*) &&\coloneqq \mathbf{N}_{\mathbf{F}}^{\textnormal{imex}}(T^*,I,E,\mathbf{F}),\\
    \frac{\partial T}{\partial t} &= - \frac{1}{\rho c_v}\left(ac\sigma_P(T^*)T^4 - c\sigma_E(T^*,I)E\right) &&\coloneqq \mathbf{N}_T^{\textnormal{imex}}(T^*,I,E,T).
\end{align}
\end{subequations}
Note that now we have removed the implicit dependence of $I$ on temperature. As a result, each implicit stage requires one transport sweep in angular intensity for fixed temperature, followed by closing and solving the (gray) nonlinear LO system in $\{E\mathbf{F},T\}$ with collapsed opacities. Collapsed opacities here are fixed and no longer depend on $T$, only on $T^*$. This is because in a single LO iteration, opacities are not updated, and by breaking the implicit coupling between $I$ and the LO system, we no longer have an outer iteration in which to converge collapsed opacities as in the semi-implicit setting with partitioned collapsed opacities.

\section{Adaptive Runge Kutta time integration}\label{sec:time}

\subsection{Adaptivity}\label{sec:time:adap}

The local error produced from one step of \eqref{eq:limex} can be estimated with the help of an embedded method of lower order. To output such an embedding, we add an additional $\widehat{\mathbf{b}}$ vector to each of the tableaux in \eqref{eq:imex-tableaux}. Like the primary method, we utilize the same embedded vector for the implicit and explicit schemes. The embedded solution is then computed analogous to the solution update in \eqref{eq:limex-sol} via
\begin{equation*}
    \widehat{y}_{n+1} = y_n + \Delta t \sum_{j=1}^s {\widehat{b}_{j}} \mathcal{N}(Y^*_{(j)},Y_{(j)}).
\end{equation*}
We denote the embedded order of accuracy with $\widehat{p}$ and the primary order of accuracy with $p$.
Following the standard approaches of \cite{Ceschino.1961} and \cite[pp. 167-168]{Hairer.1993}, a new timestep is selected with
\begin{equation*}
    \Delta t_{new} = 0.9 \Delta t \cdot err^{-1 / (\min(p, \widehat{p}) + 1)},
    \qquad
    err = \sqrt{\frac{1}{N} \sum_{i = 1}^N \left( \frac{y_{n+1,i} - \widehat{y}_{n+1,i}}{atol_i + rtol_i \max(|y_{n,i}|, |y_{n+1,i}|)} \right)^2},
\end{equation*}
and a step is deemed acceptable for the specified absolute and relative tolerances if $err \leq 1$. Note, the error estimate is computed with a local relative tolerance to ensure that error associated with small solution values is not ignored, something that is particularly important in TRT problems where solutions can span many orders of magnitudes, and accurate wavefront evolution is critical to long-term accuracy of dynamics. 

In this work, we focus on methods with order $p = 2$ and an embedding one order lower ($\widehat{p} = 1$). We use second-order methods to limit the memory requirements, and because our experience with TRT and radiation hydrodynamics has indicated no benefit in terms of accuracy vs. wallclock time of higher order methods. To the authors' knowledge, IMEX embedded pairs of order two have only been constructed in \cite{Giraldo.2013}, although higher order methods can be found in \cite{Kennedy.2003tv4,Kennedy.2019}. 
We are particularly interested in schemes that are asymptotically preserving for hyperbolic equations with stiff relaxation terms in the sense of \cite{Pareschi.2005}, a property that our previous work has found important for TRT \cite{imex-trt} and radiation hydrodynamics \cite{rad-hydro}. We have also found EDIRK schemes as in \cite{Kennedy.2003tv4,Kennedy.2019} to be numerically unstable for the very stiff physics in TRT. There are few schemes that satisfy these property and have derived embeddings in the literature, and to the best of our knowledge, we are aware of no such second order schemes (as we are interested in) with derived embeddings. With the sparse availability of existing embedded pairs with the desired properties, we look to add embeddings to existing methods from the literature. For a given primary method, an embedding is typically not unique. Free parameters can be used to optimize stability and a number of measures for embedding quality. Before defining these measures, we must first present the structure of the local truncation error of the semi-implicit RK scheme \eqref{eq:limex}. This can be derived via Taylor series, or more systematically with the tree-based approach of \cite{nprk2}:
\begin{equation} \label{eq:local_error}
    \begin{split}
        y(t_1) - y_1
        &= \Delta t \tau^{(1)}_1 \mathcal{N}_0
        + \Delta t^2 \tau^{(2)}_1 D_1 \mathcal{N}_0 \cdot \mathcal{N}_0
        + \Delta t^2 \tau^{(2)}_2 D_2 \mathcal{N}_0 \cdot \mathcal{N}_0
        + \Delta t^3 \tau^{(3)}_1 D_1 \mathcal{N}_0 \cdot D_1 \mathcal{N}_0 \cdot \mathcal{N}_0 \\
        & \quad + \Delta t^3 \tau^{(3)}_2 D_1 \mathcal{N}_0 \cdot D_2 \mathcal{N}_0 \cdot \mathcal{N}_0
        + \Delta t^3 \tau^{(3)}_3 D_2 \mathcal{N}_0 \cdot D_1 \mathcal{N}_0 \cdot \mathcal{N}_0
        + \Delta t^3 \tau^{(3)}_4 D_2 \mathcal{N}_0 \cdot D_2 \mathcal{N}_0 \cdot \mathcal{N}_0 \\
        & \quad + \Delta t^3 \tau^{(3)}_5 D_{11} \mathcal{N}_0(\mathcal{N}_0, \mathcal{N}_0)
        + \Delta t^3 \tau^{(3)}_6 D_{12} \mathcal{N}_0(\mathcal{N}_0, \mathcal{N}_0)
        + \Delta t^3 \tau^{(3)}_7 D_{22} \mathcal{N}_0(\mathcal{N}_0, \mathcal{N}_0)
        + \mathcal{O}(\Delta t^4)
    \end{split}
\end{equation}
For brevity, we have used a subscript of $0$ to denote evaluation at the initial condition, e.g., $\mathcal{N}_0 = \mathcal{N}(y_0, y_0)$. Differentiation of $\mathcal{N}$ with respect to the $i$-th argument is denoted by $D_i$. Finally, $\tau_i^{(j)}$ represents the $i$-th order condition residual of order $j$ with values given in Table \ref{tab:order_conditions}. Analogously, residuals for the embedded method use a hat, e.g., $\widehat{\tau}_1^{(1)} = 1 - \widehat{\mathbf{b}}^T \mathbf{e}$ where $\mathbf{e} = [1, \dots, 1]^T \in \mathbb{R}^{s}$.

\begin{table}[ht!]
    \centering
    \begin{tabular}{l|l|lll}
        Order 1 & Order 2 & Order 3 \\ \hline
        $\tau^{(1)}_1 = 1 - \mathbf{b}^T \mathbf{e}$
        & $\tau^{(2)}_1 = \frac{1}{2} - \mathbf{b}^T \mathbf{c}$,
        & $\tau^{(3)}_1 = \frac{1}{6} - \mathbf{b}^T A \mathbf{c}$,
        & $\tau^{(3)}_2 = \frac{1}{6} - \mathbf{b}^T A \tilde{\mathbf{c}}$,
        & $\tau^{(3)}_3 = \frac{1}{6} - \mathbf{b}^T \tilde{A} \mathbf{c}$, \\
        & $\tau^{(2)}_2 = \frac{1}{2} - \mathbf{b}^T \tilde{\mathbf{c}}$
        & $\tau^{(3)}_4 = \frac{1}{6} - \mathbf{b}^T \tilde{A} \tilde{\mathbf{c}}$,
        & $\tau^{(3)}_5 = \frac{1}{6} - \frac{1}{2} \mathbf{b}^T \mathbf{c}^2$,
        & $\tau^{(3)}_6 = \frac{1}{3} - \mathbf{b}^T (\mathbf{c} \times \tilde{\mathbf{c}})$, \\
        & & $\tau^{(3)}_7 = \frac{1}{6} - \frac{1}{2} \mathbf{b}^T \tilde{\mathbf{c}}^2$
    \end{tabular}
    \caption{Order condition residuals up to order three for the local error expansion \eqref{eq:local_error}.}
    \label{tab:order_conditions}
\end{table}

Following \cite{Kennedy.2003tv4}, the principal error for the primary an embedded method are
\begin{equation*}
    A^{(p+1)} = \sqrt{\sum_{i} \left( \tau^{(p+1)}_i \right)^2},
    \quad
    \text{and}
    \quad
    \widehat{A}^{(\widehat{p}+1)} = \sqrt{\sum_{i} \left( \widehat{\tau}^{(\widehat{p}+1)}_i \right)^2},
\end{equation*}
respectively. To ensure the embedded method provides an accurate error estimate over a wide range of $\Delta t$, we aim to have the following quantities close to one:
\begin{equation*}
    B^{(\widehat{p}+2)} = \frac{\widehat{A}^{(\widehat{p}+2)}}{\widehat{A}^{(\widehat{p}+1)}},
    \quad
    C^{(\widehat{p}+2)} = \frac{\sqrt{\sum_{i} \left( \widehat{\tau}^{(\widehat{p}+2)}_i - \tau^{(\widehat{p}+2)}_i \right)^2}}{\widehat{A}^{(\widehat{p}+1)}},
    \quad
    E^{(\widehat{p}+2)} = \frac{A^{(\widehat{p}+2)}}{\widehat{A}^{(\widehat{p}+1)}}.
\end{equation*}
Analogous measures of embedding quality for the underlying implicit and explicit Runge--Kutta methods are denoted with a an $I$ and $E$ subscript, respectively.

We formulate the task of computing an embedded method as follows:
\begin{equation} \label{eq:optimal_embedding}
    \begin{aligned}
        \argmin_{\widehat{b}} \quad &
        \left\| [B^{(3)}, C^{(3)}, E^{(3)}, B^{(3)}_I, C^{(3)}_I, E^{(3)}_I, B^{(3)}_E, C^{(3)}_E, E^{(3)}_E]^T - [1, \dots, 1]^T \right\|_2 \\
        \textrm{s.t.} \quad & 0 = \widehat{\tau}_1^{(1)}
        \text{ and embedding of underlying implicit method is A-stable} \\
    \end{aligned}
\end{equation}
\tcb{A stronger constraint of L-stability in \eqref{eq:optimal_embedding} would be advantageous, but unfortunately, it is infeasible for the methods discussed in the following section.}
We use Mathematica and the Integreat library \cite{Roberts.2023} to compute the measures of embedding quality, an algebraic expression for the A-stability constraint, and a global minimizer numerically.

\subsection{Embedded IMEX-RK schemes}\label{sec:time:schemes}

Following our recent papers \cite{imex-trt,rad-hydro}, we have found four IMEX-RK tableaux pairs to be significantly more robust than all others we have tested for the very stiff and nonlinear problems that arise in TRT and its coupling to other physics. \tcb{Joint stability for (nonlinear) partitioned integration of nonlinear equations is largely unexplored, and what precisely makes for a robust scheme remains an open question. However, we have found the asymptotic preserving property for hyperbolic systems with relaxation (a model that can be seen as a natural simplification of TRT and radiation-hydrodynamics) as defined in \cite{Pareschi.2005} to be consistently useful/important. In contrast, although EDIRK methods \cite{Kennedy.2003tv4,Kennedy.2019} can offer higher stage order, we have not found them to be robust for problems in TRT or rad-hydro. We believe this is due to instabilities arising from the first explicit step over highly stiff radiation physics that are not adequately damped in subsequent stages (part of this may be numerical in nature due to the large scale variations). We have also generally found L-stable implicit schemes and favorable joint linear stability properties as we analyze in \cite{nprk1,buvoli2025multirate} to be important. For example, IMEX schemes based on an implicit midpoint method have tended to have poor stability in practice, likely due to an only A-stable implicit component and poor joint linear stability.}

All schemes we utilize are second order with L-stable implicit components, and three of the four were designed specifically for hyperbolic systems with relaxation \cite{Pareschi.2005}. The fourth is an equivalent additive representation of a certain scheme in the new class of nonlinearly partitioned Runge--Kutta methods (NPRK) \cite{nprk1,nprk2}. Here we present each scheme, including an embedding which is a nearby rational approximation to the optimal embedding \eqref{eq:optimal_embedding}.
\tcb{Note that the embedded methods inherit the asymptotic preserving property since $A$ is nonsingular for all methods (see \cite[Theorem 3.1]{Pareschi.2005}).}
We summarize some of the other method properties in Table \ref{tab:IMEX_properties} \tcb{and provide stability plots in \ref{app:stability}}.
\newline\newline\noindent\textbf{H-LDIRK2(2,2,2)} \cite[Table II]{Pareschi.2005}:
\begin{align*}
    \begin{array}{c | c c c}
    0 & 0 & 0 \\
    1 & 1 & 0 \\\hline
    & 1/2 & 1/2 \\\hline
    & 3/10 & 7/10
    \end{array},
\hspace{5ex}
    \begin{array}{c | c c}
    \gamma & \gamma & 0 \\
    1 -\gamma & 1-2\gamma & \gamma \\\hline
    & 1/2 & 1/2 \\\hline
    & 3/10 & 7/10
    \end{array},
\hspace{2ex} \textnormal{for }\gamma=1-1/\sqrt{2}.
\end{align*}
\noindent\textbf{SSP-LDIRK2(3,3,2)} \cite[Table IV]{Pareschi.2005}:
\begin{align*}
    \begin{array}{c | c c c}
    0 & 0 & 0 & 0 \\
    1/2 & 1/2 & 0 & 0\\
    1 & 1/2 & 1/2 & 0 \\\hline
    & 1/3 & 1/3 & 1/3 \\\hline
    & 7/41 & 13/33 & 589/1353
    \end{array},
\hspace{5ex}
    \begin{array}{c | c c c}
    1/4 & 1/4 & 0 & 0 \\
    1/4 & 0 & 1/4 & 0 \\
    1 & 1/3 & 1/3 & 1/3 \\\hline
    & 1/3 & 1/3 & 1/3 \\\hline
    & 7/41 & 13/33 & 589/1353
    \end{array}.
\end{align*}
\noindent\textbf{SSP-LDIRK3(3,3,2)}\cite[Table V]{Pareschi.2005}:
\begin{align*}
    \begin{array}{c | c c c}
    0 & 0 & 0 & 0 \\
    1 & 1 & 0 & 0\\
    1/2 & 1/4 & 1/4 & 0 \\\hline
    & 1/6 & 1/6 & 2/3 \\\hline
    & -2/17 & 2/13 & 213/221
    \end{array},
\hspace{5ex}
    \begin{array}{c | c c c}
    \gamma & \gamma & 0 & 0 \\
    1-\gamma & 1-2\gamma & \gamma & 0 \\
    1/2 & 1/2 - \gamma & 0 & \gamma \\\hline
    & 1/6 & 1/6 & 2/3 \\\hline
    & -2/17 & 2/13 & 213/221
    \end{array},
    \hspace{2ex} \textnormal{for }\gamma=1-1/\sqrt{2}.
\end{align*}

\noindent\textbf{IMEX-NPRK2[42]b} \cite[Eq. SM5.35]{nprk1}:
\begin{align*}
    \begin{array}{c | c c c}
     & 0 & 0 \\
     & \frac{26+3\sqrt{2}}{42} & 0 \\\hline
    & \frac{16 + 9\sqrt{2}}{94} & \frac{78 - 9\sqrt{2}}{94}\\\hline
    & 1/15 & 14/15
    \end{array},
\hspace{5ex}
    \begin{array}{c | c c}
     & \gamma & 0 \\
     & \frac{-20 + 23\sqrt{2}}{42} & \gamma \\\hline
    & \frac{16 + 9\sqrt{2}}{94} & \frac{78 - 9\sqrt{2}}{94} \\\hline
    & 1/15 & 14/15
    \end{array},
\hspace{2ex} \textnormal{for }\gamma=1-1/\sqrt{2}.
\end{align*}

\begin{table}[ht!]
    \centering
    \begin{tabular}{r|ccccccc}
        Method & $s$ & \shortstack{Implicit Method\\Stability} & \shortstack{Implicit Embedded\\Stability} & $A^{(3)}$ & $B^{(3)}$ & $C^{(3)}$ & $E^{(3)}$ \\
        \hline
        H-LDIRK2(2,2,2) & 2 & L-stable & A-stable & 0.2757 & 1.452 & 0.9346 & 0.9004 \\
        SSP-LDIRK2(3,3,2) & 3 & L-stable & A-stable & 0.1792 & 1.583 & 1.032 & 0.8292 \\
        SSP-LDIRK3(3,3,2) & 3 & L-stable & A-stable & 0.1265 & 1.116 & 0.9339 & 0.6088 \\
        IMEX-NPRK2[42]b & 2 & L-stable & A-stable & 0.2564 & 1.153 & 0.7964 & 0.9736 \\
    \end{tabular}
    \caption{Properties of the second order IMEX-RK schemes and their embeddings. A smaller value of $A^{(3)}$ typically indicates a more accurate method. Ideally, $B^{(3)}$, $C^{(3)}$, $E^{(3)}$ should be close to one.}
    \label{tab:IMEX_properties}
\end{table}

\subsection{Moment-based adaptivity for TRT}\label{sec:time:moment}

Due to the high dimensionality of the TRT solution space, memory is a major concern in numerical simulation. As a result, adding an additional full copy of the specific intensity for purposes of error estimation (in addition to an extra copy or two for a second-order integrator) is not desirable. To that end, we propose performing local error estimation only in the reduced moment system, particularly in the radiation energy density and the temperature. In addition to being relatively cheap to store an additional copy for error estimation, these variables are also the primary coupling components to other physics, such as in radiation hydrodynamics. Thus, we believe these lower dimensional variables are both computationally feasible and physically suitable representations of error dynamics to guide time step selection. In the error estimation algorithm from \Cref{sec:time:adap}, we thus restrict our embedded vector to only consider the temperature and/or radiation energy components. 

Note, one could also consider using radiation flux $\mathbf{F}$ for error estimation. We have found $\mathbf{F}$ to consistently drive timesteps many orders of magnitude smaller than needed. We believe this is due to inconsistency in the flux closure compared with the flux integrated directly from $I$, either from the IMEX time integration or inexact solution of the semi-implicit system. This is consistent with results seen for linear neutron transport in \cite[Fig. 5.8]{olivier2024consistent}, where inconsistent moment discretizations can lead to nontrivial error in the first moment.

\section{Numerical results}\label{sec:results}

We demonstrate our moment-based adaptive time-stepping schemes using the newly derived embedded IMEX pairs on two common test problems in TRT: \tcb{(i) the gray tophat problem to test 2d multi-material geometries, and (ii) the multifrequency Larsen problem to test multifrequency TRT in highly non-equilibrium regimes.}

\subsection{The 2d Tophat Problem}\label{sec:results:tophat}

\begin{figure}[!htbp]
    \centering
    \begin{minipage}[b]{0.2\textwidth}
        \centering
        \includegraphics[width=\textwidth]{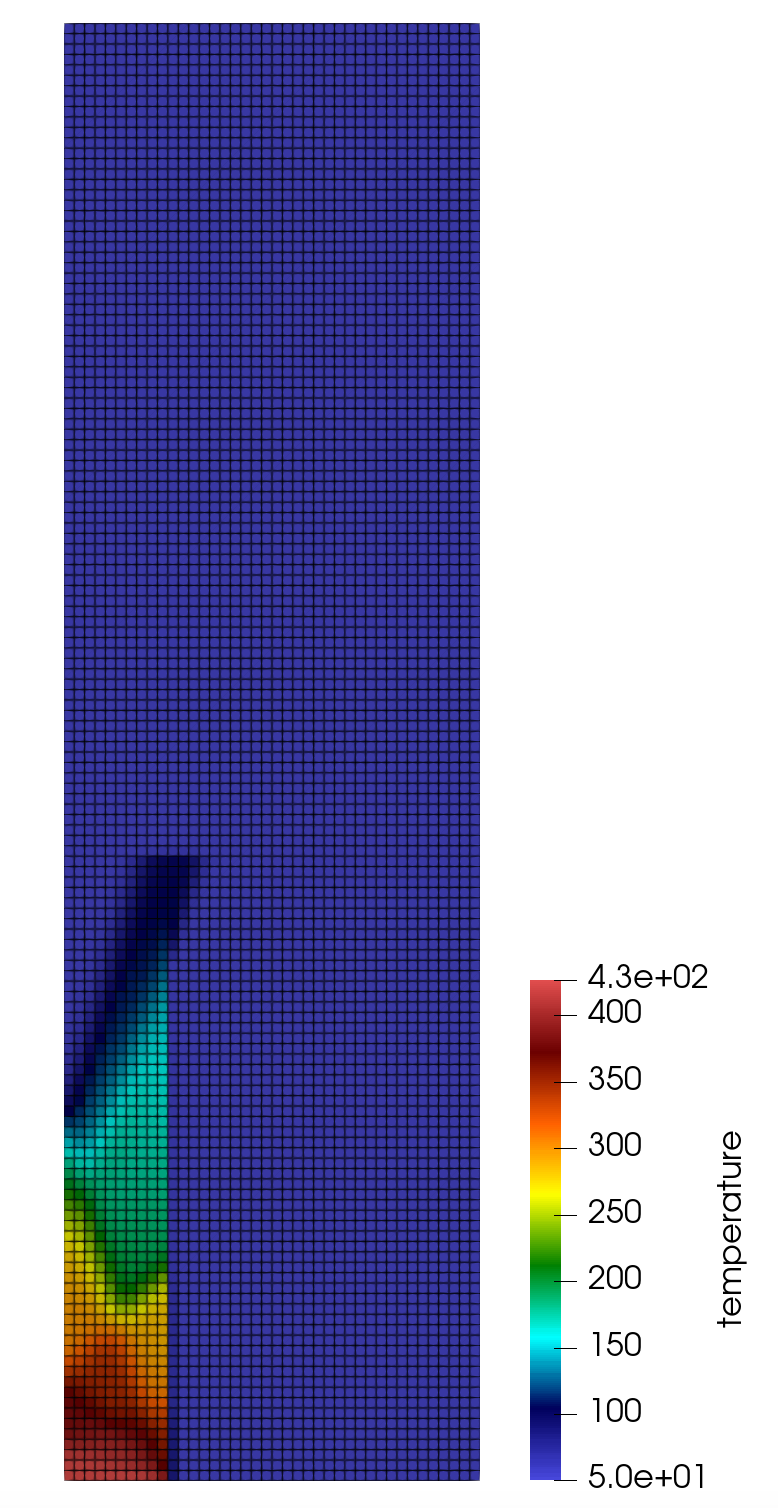} 
        \caption*{(a)}
    \end{minipage}
    \begin{minipage}[b]{0.2\textwidth}
        \centering
        \includegraphics[width=\textwidth]{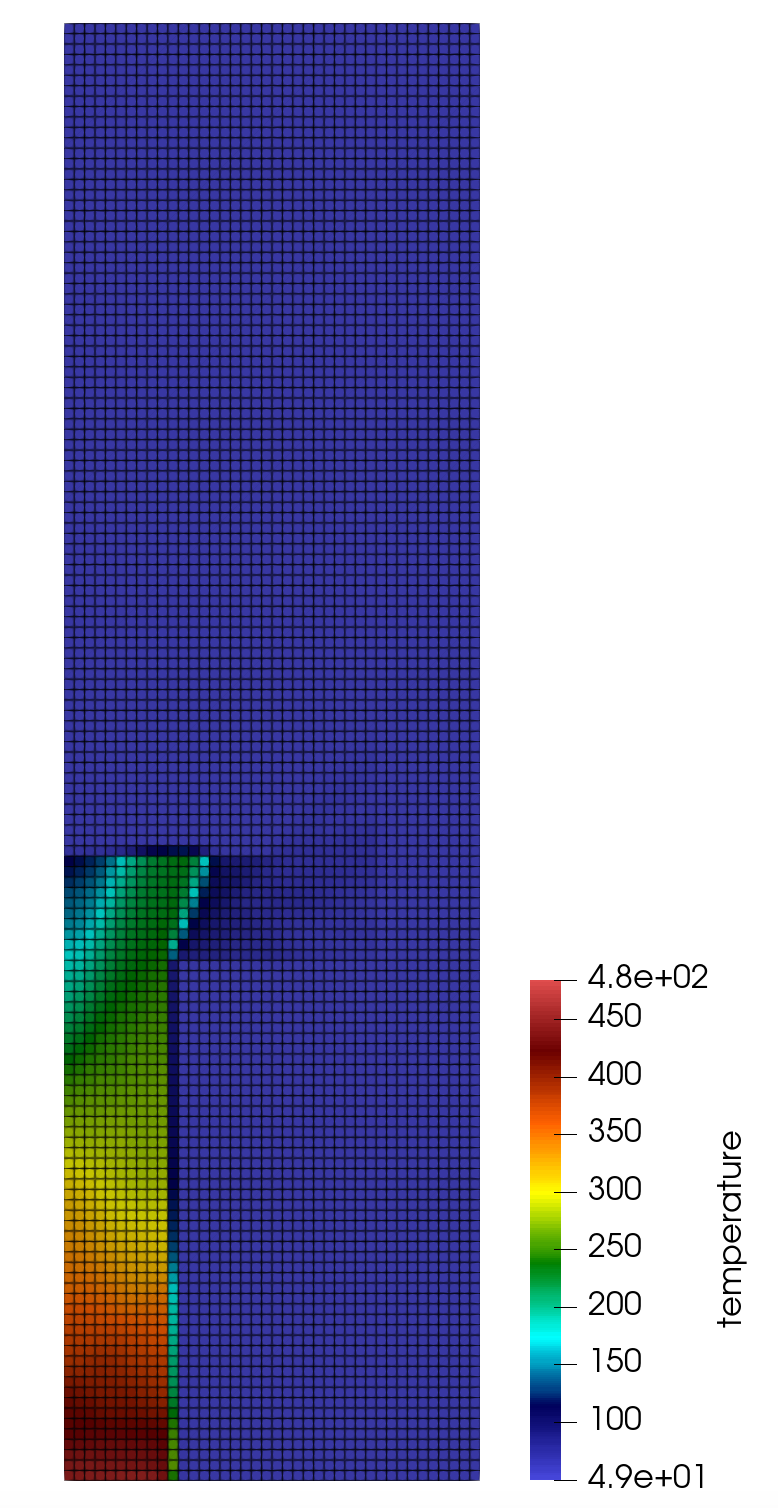} 
        \caption*{(b)}
    \end{minipage}
    \begin{minipage}[b]{0.2\textwidth}
         \centering
        \includegraphics[width=\textwidth]{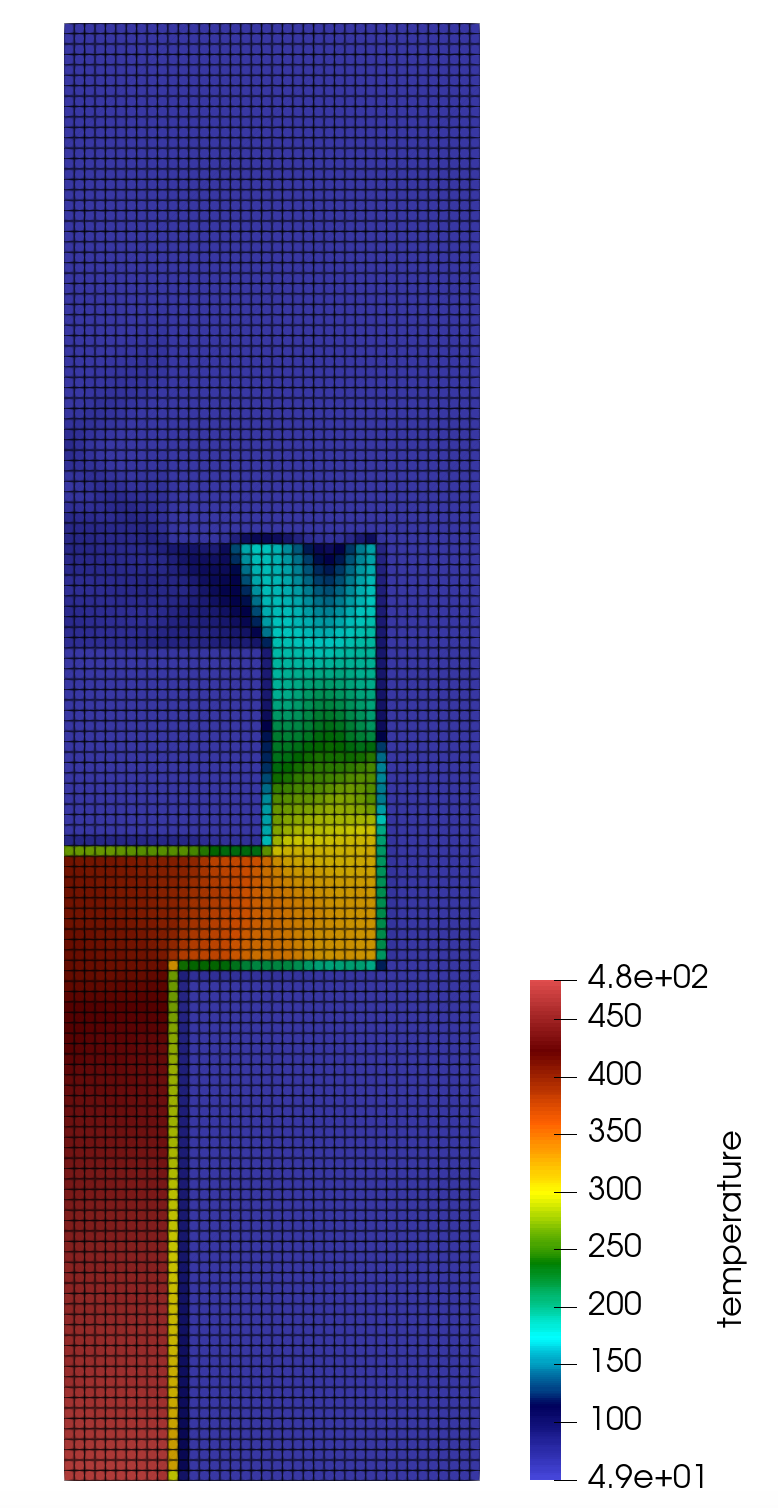} 
        \caption*{(c)}
    \end{minipage}
    \begin{minipage}[b]{0.2\textwidth}
        \centering
        \includegraphics[width=\textwidth]{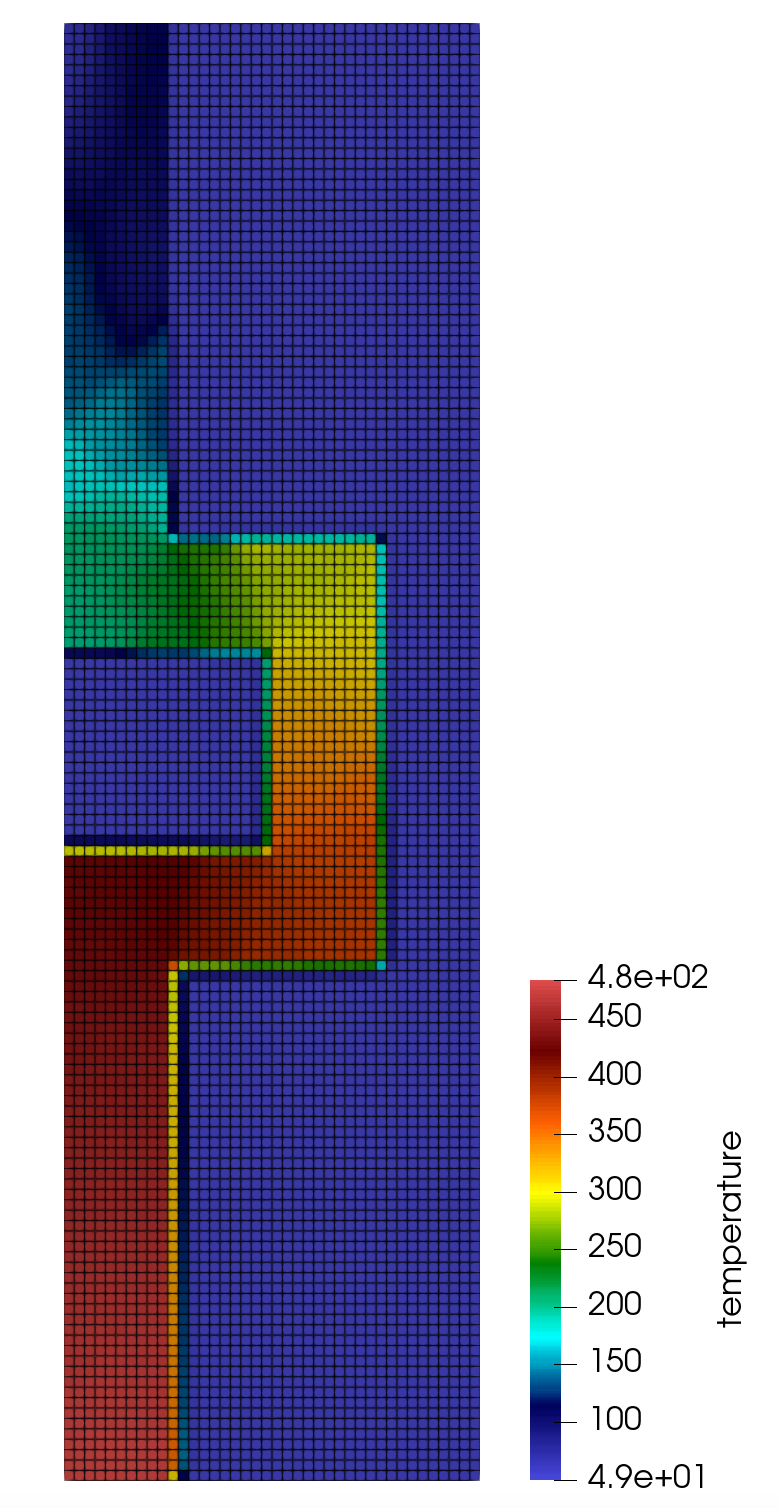} 
        \caption*{(d)}
    \end{minipage}

    \caption{Evolution of the temperature field in the 2d Tophat problem.  The times shown are (a) $t=5e-10$, (b) $t=5e-09$, (c) $t=5e-08$ and (d) $t=1e-07$. }
    \label{fig:temperature_field_tophat_IMEX103}
\end{figure}

First we consider the multimaterial tophat problem. The problem consists of two materials: the optically thin pipe and optically thick wall described by $\sigma_\text{pipe} = \SI{0.2}{\per\cm}$ and $\sigma_\text{wall} = \SI{2000}{\per\cm}$, respectively, which do not depend on temperature and are thus fixed in time. The heat capacity is $C_v = 10^{12}\si{\erg\per\eV\per\cm\cubed}$ in both the wall and pipe. Radiation enters the pipe at the bottom of the domain according to a Planckian distribution evaluated at a temperature of $T_b= $ \SI{500}{\eV}, that is $I_\text{inflow} = \frac{ac T_b^4}{4\pi}$. Elsewhere, the domain boundary is treated as a vacuum. The computational domain is halved by applying a reflection boundary along the plane $x=0$. The initial temperature and radiation fields are set to be in equilibrium with each other at a spatially uniform temperature of \SI{50}{\eV}. The evolution of the temperature field through the pipe can be seen in \Cref{fig:temperature_field_tophat_IMEX103}. We run this problem on a uniform grid of $40\times140$ cells with 6th-order S$_N$ quadrature and HOLO discretization as in \cite{Park.2020,imex-trt}.

\begin{figure}[!htbp]
    \centering
    \begin{minipage}{0.75\textwidth} 
        \centering
        \begin{minipage}[b]{0.48\textwidth} 
            \centering
            \includegraphics[width=\textwidth]{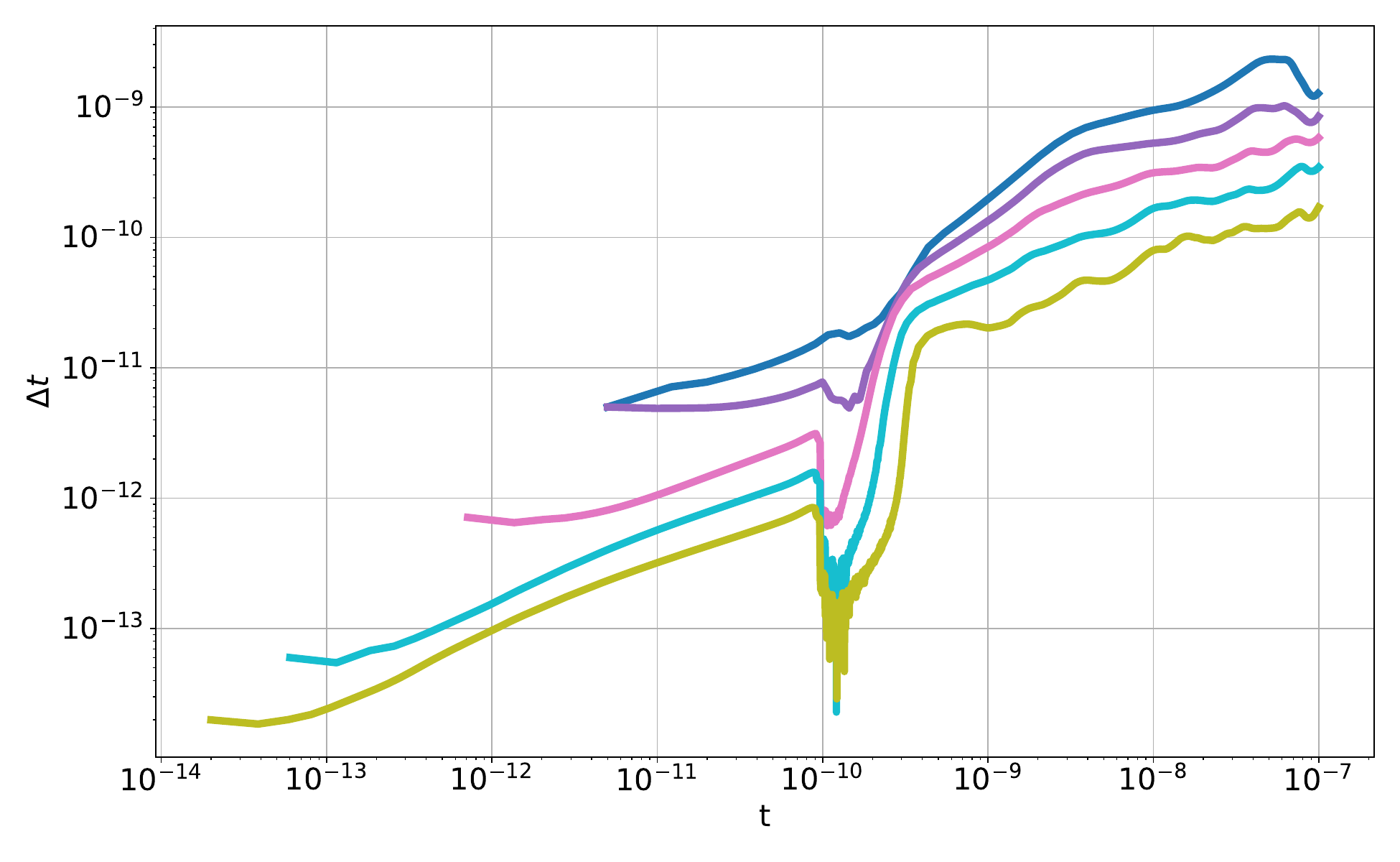}
            \caption*{(a) H-LDIRK2(2,2,2)}
        \end{minipage}
        \hfill
        \begin{minipage}[b]{0.48\textwidth} 
            \centering
            \includegraphics[width=\textwidth]{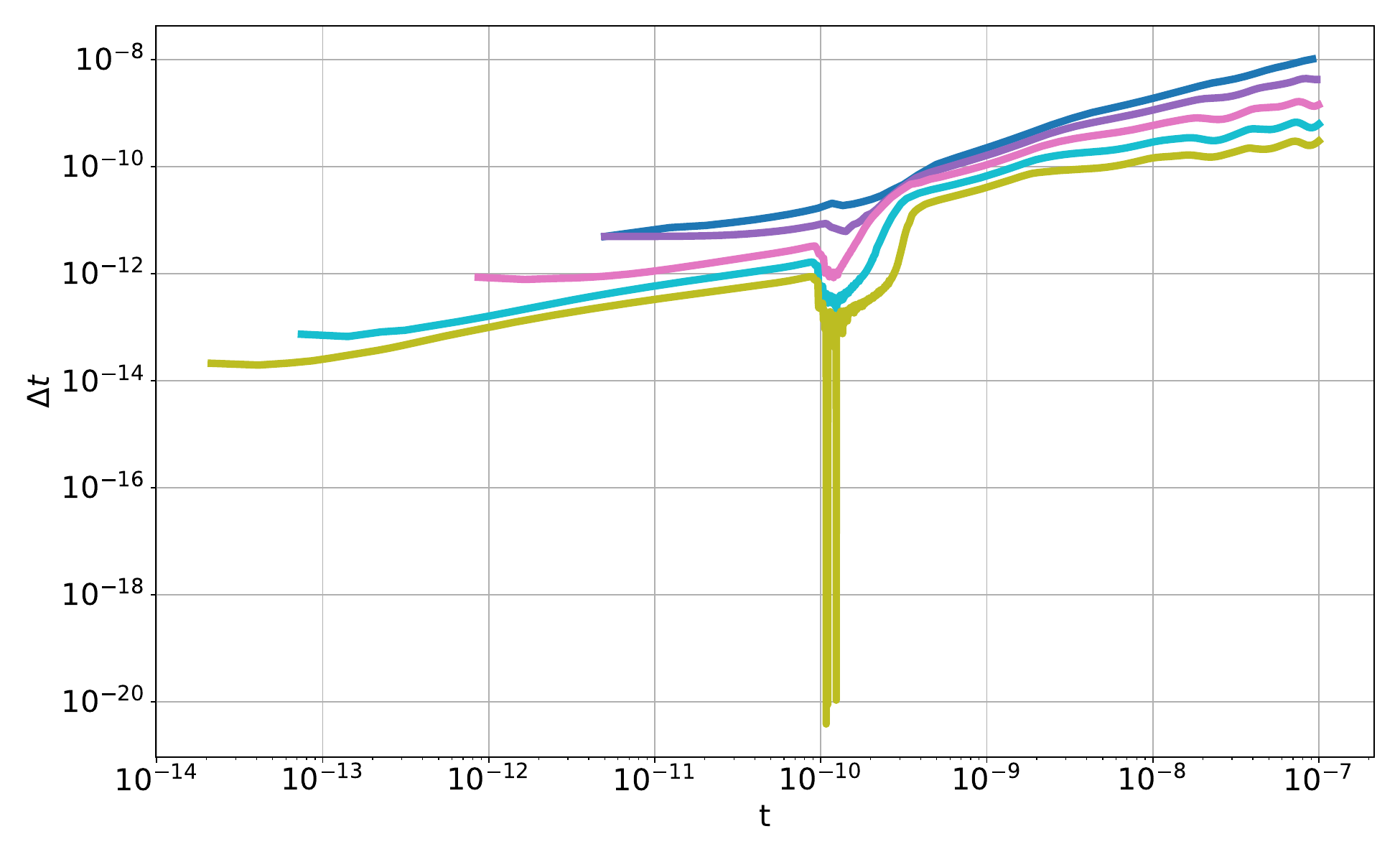}
            \caption*{(b) SSP-LDIRK2(3,3,2)}
        \end{minipage}

        \vspace{5pt} 

        \begin{minipage}[b]{0.48\textwidth} 
            \centering
            \includegraphics[width=\textwidth]{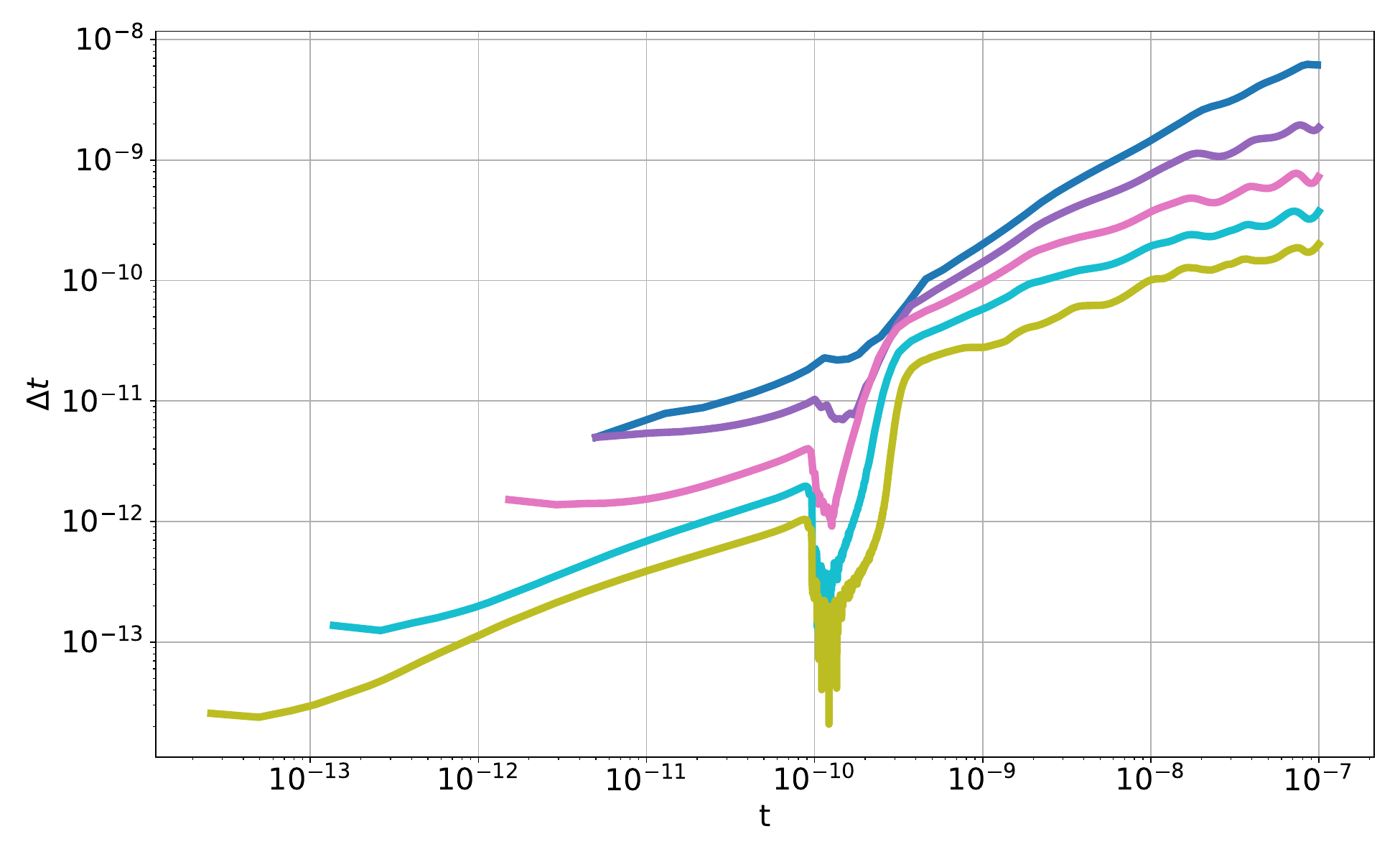}
            \caption*{(c) SSP-LDIRK3(3,3,2)}
        \end{minipage}
        \hfill
        \begin{minipage}[b]{0.48\textwidth} 
            \centering
            \includegraphics[width=\textwidth]{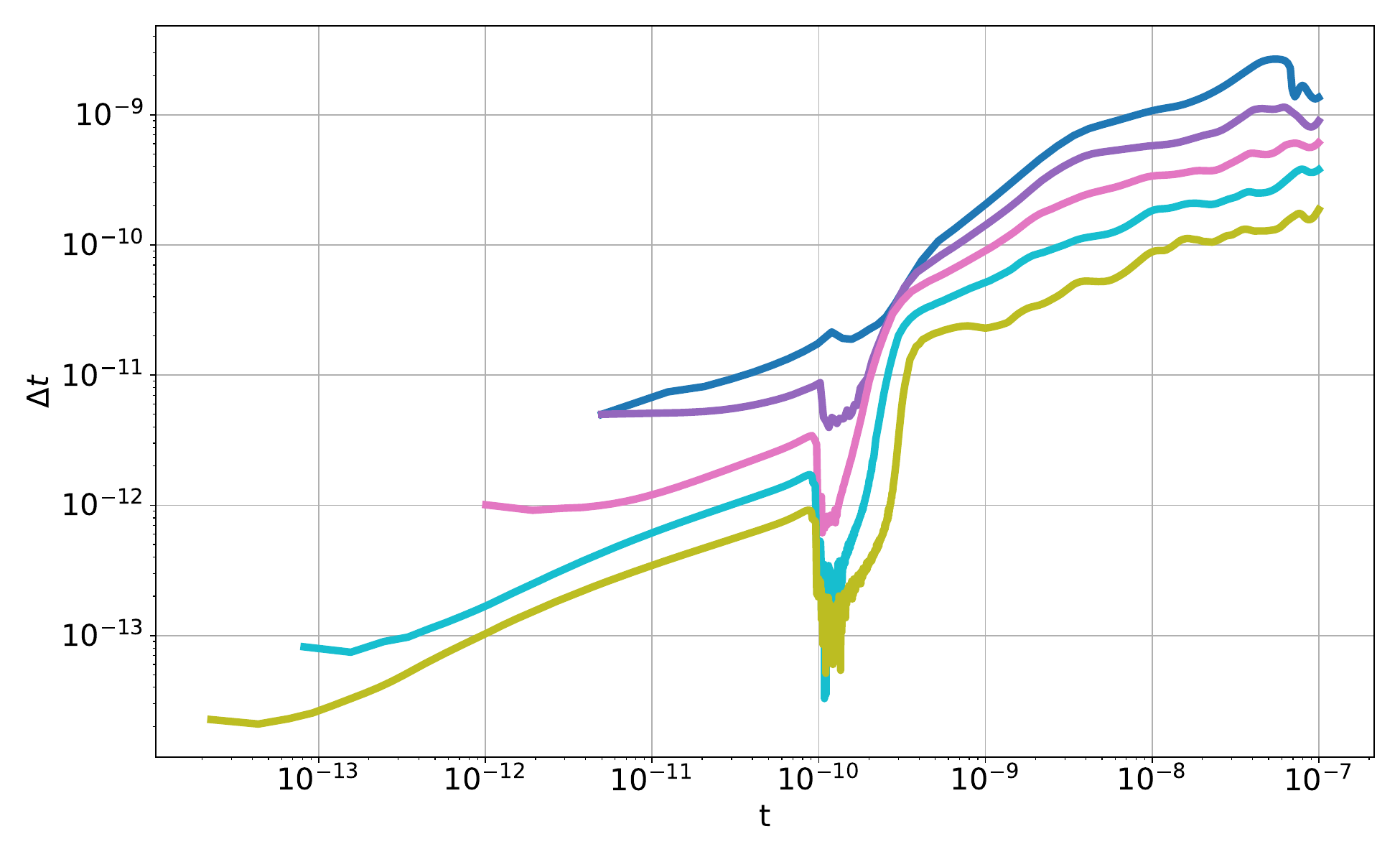}
            \caption*{(d) IMEX-NPRK[42]b}
        \end{minipage}
    \end{minipage}
    \begin{minipage}{0.18\textwidth} 
        \centering
        \includegraphics[width=\textwidth]{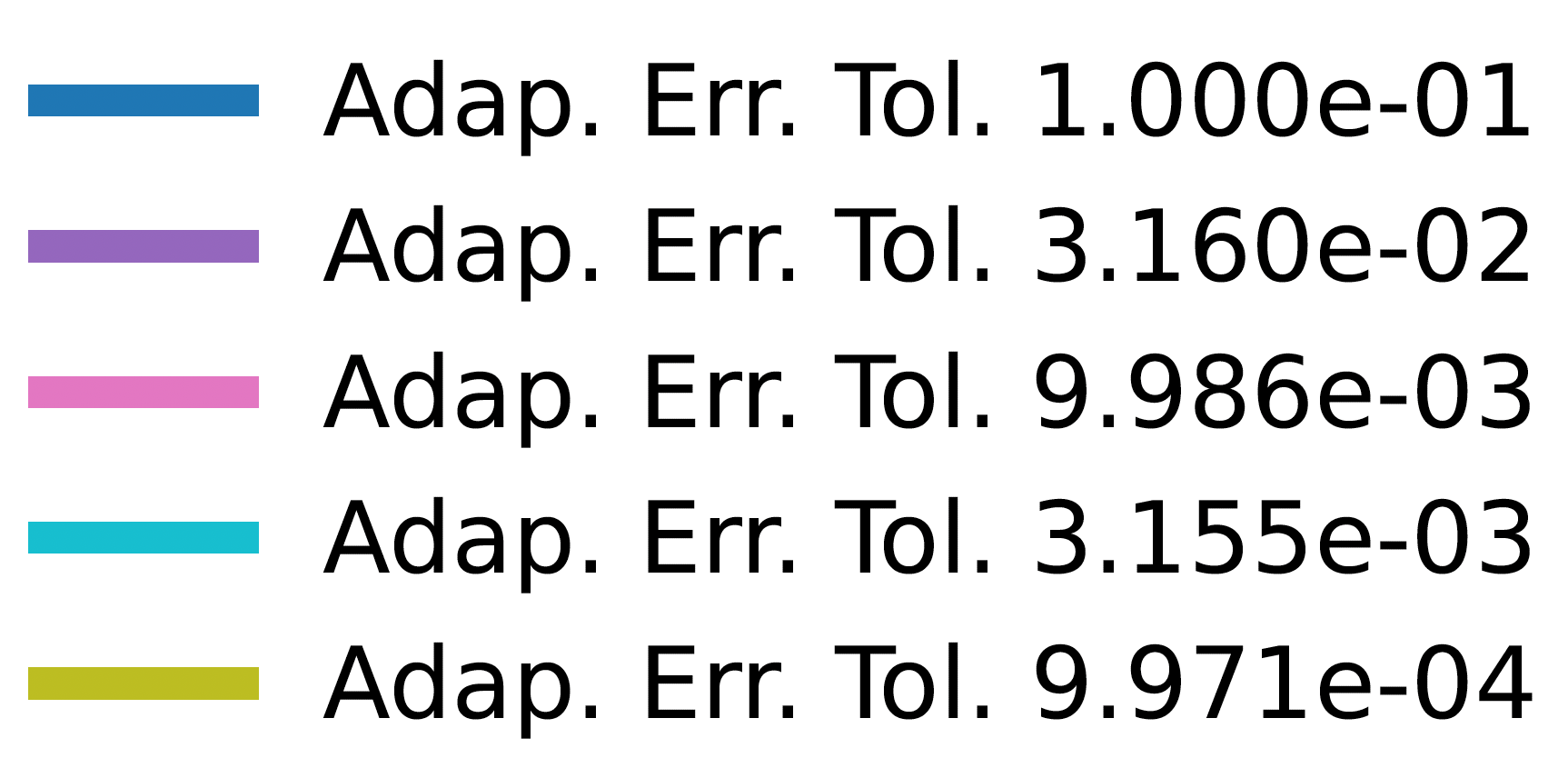}\\\includegraphics[height=0.3\textheight]{Er_err_paraview_output_tophat_tophat_imex103_temp_t5e-10.png}
        \caption*{(e)}
    \end{minipage}

    \caption{Plot of $\Delta t$ vs. $t$ for IMEX schemes for the 2d Tophat problem (a) H-LDIRK2(2,2,2), (b) SSP-LDIRK2(3,3,2) (c) SSP-LDIRK3(3,3,2) and (d) IMEX-NPRK[42]b.  In (e), the temperature field is shown at $t=5e-10$ as the field has to resolve the first corner.  The resolution of the first corner corresponds with the drop in $\Delta t$ for each of the considered schemes. }
    \label{fig:IMEX_dt_v_t_plots_with_temp_at_t5e-10}
\end{figure}

We first demonstrate that the adaptive time-stepping naturally captures the dynamic timescale of the problem, and is  robust across time integration scheme and formulation (IMEX vs semi-implicit). \Cref{fig:IMEX_dt_v_t_plots_with_temp_at_t5e-10} shows the adaptive time step choice as a function of physical time in simulation in seconds for IMEX integration, for the four IMEX schemes and embedded error estimators, and four different adaptive error tolerances. Error here is measured with respect to radiation energy density. The adaptively chosen timestep is intuitive as well; early time when radiation is entering the domain and the problem is streaming dominated, the timestep must be very small to capture the evolution of the radiation energy and temperature front. As the pipe heats up, the timestep begins to increase. Around physical time $t=10^{-10}$s, the temperature starts turning the first corner in the pipe (see \Cref{fig:IMEX_dt_v_t_plots_with_temp_at_t5e-10}e). Here, all schemes lead to a drop in timestep to capture these dynamics, with a more significant decrease in timestep required for higher accuracy. After the corner is resolved, the full pipe begins to heat up. Eventually, the pipe becomes saturated and dynamics transition to a slow diffusive timescale, which can be captured with large time steps well over the streaming timescale. Altogether, we span 3--4 orders of mangitude in timestep size over the coarse of the simulation. 

Note that across almost all of the time domain, the adaptively chosen timestep is quite similar across all four time integration schemes and embedded error estimators. The one exception is around physical time $t=10^{-10}$s where the temperature begins turning the corner. For the smallest adaptive error tolerance, SSP-LDIRK2(3,3,2) requires minuscule timesteps on the order of $10^{-20}$s. Analyzing the simulation results, this appears to be due to negativities that arise in the embedded solution that are ``fixed,'' but such fixing is inconsistent with the solution of the high-order scheme, indicating error and a need for a smaller timestep. Note that a similar problem occurs with semi-implicit integration for SSP-LDIRK2(3,3,2). We believe this is a subtle nonlinear stability defect of the embedded IMEX pair of SSP-LDIRK2(3,3,2), and for this reason we recommend the other schemes as they have proven more robust. 

\begin{figure}[!htbp]
    \centering
    \begin{minipage}{0.75\textwidth} 
        \centering
        \begin{minipage}[b]{0.48\textwidth} 
            \centering
            \includegraphics[width=\textwidth]{Er_tvdt_tophat-imex-adap-tf1.0e-07-scheme-103-trad-dt_plot.pdf}
            \caption*{(a) IMEX, $E_r$-adaptivity}
        \end{minipage}
        \hfill
        \begin{minipage}[b]{0.48\textwidth} 
            \centering
            \includegraphics[width=\textwidth]{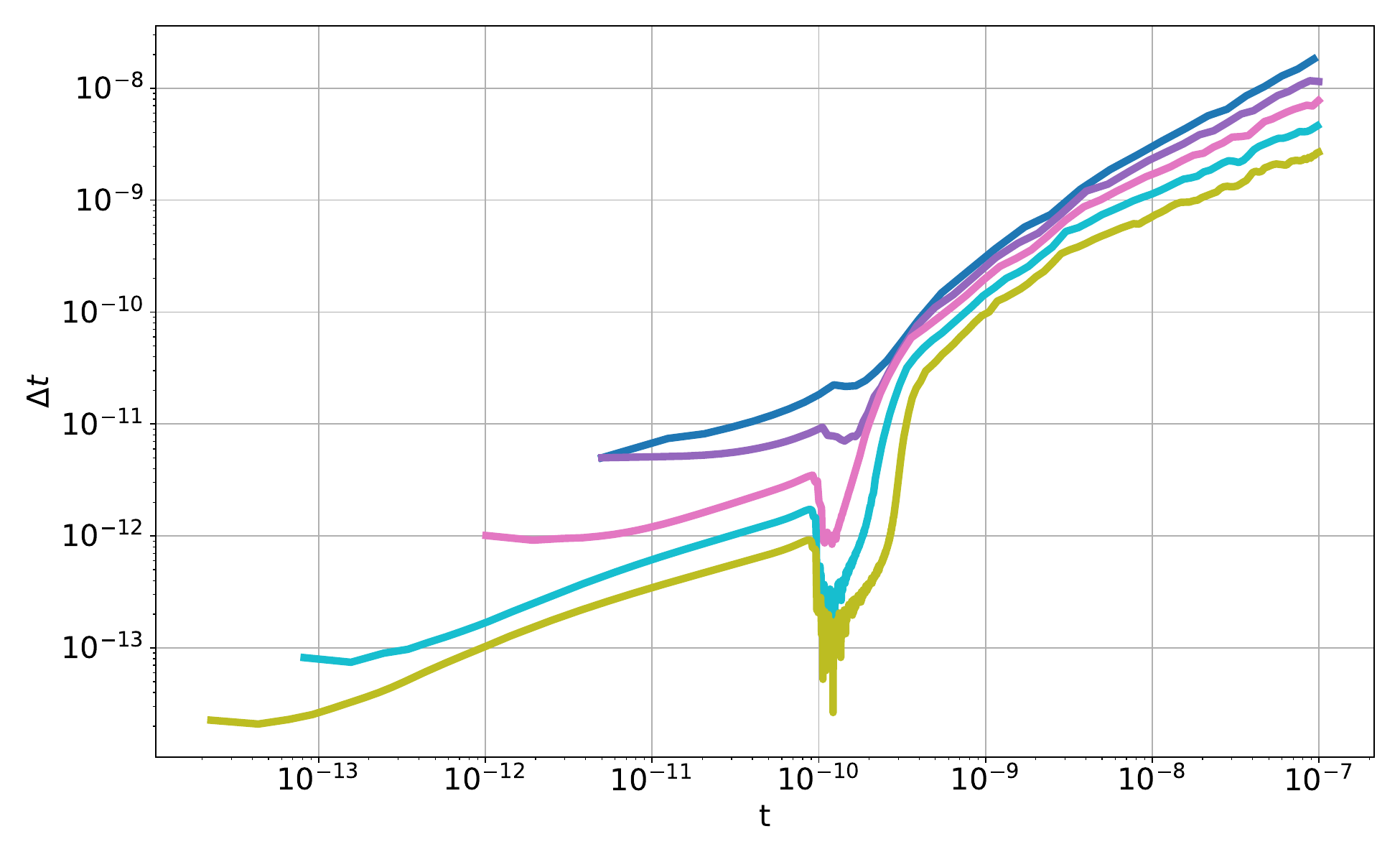}
            \caption*{(b) Semi-implicit, $E_r$-adaptivity}
        \end{minipage}

        \vspace{5pt} 

        \begin{minipage}[b]{0.48\textwidth} 
            \centering
            \includegraphics[width=\textwidth]{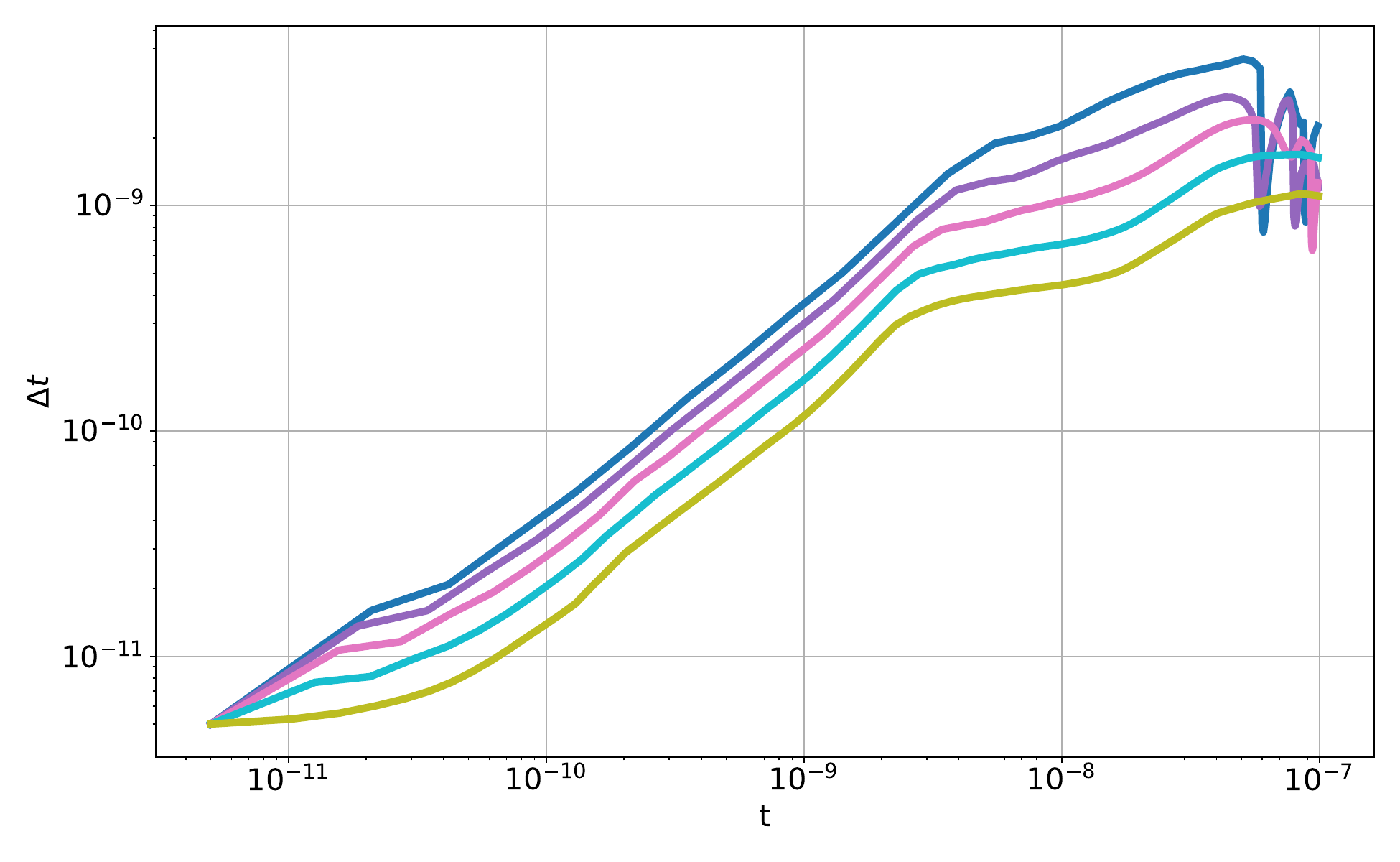}
            \caption*{(c) IMEX, $T$-adaptivity}
        \end{minipage}
        \hfill
        \begin{minipage}[b]{0.48\textwidth} 
            \centering
            \includegraphics[width=\textwidth]{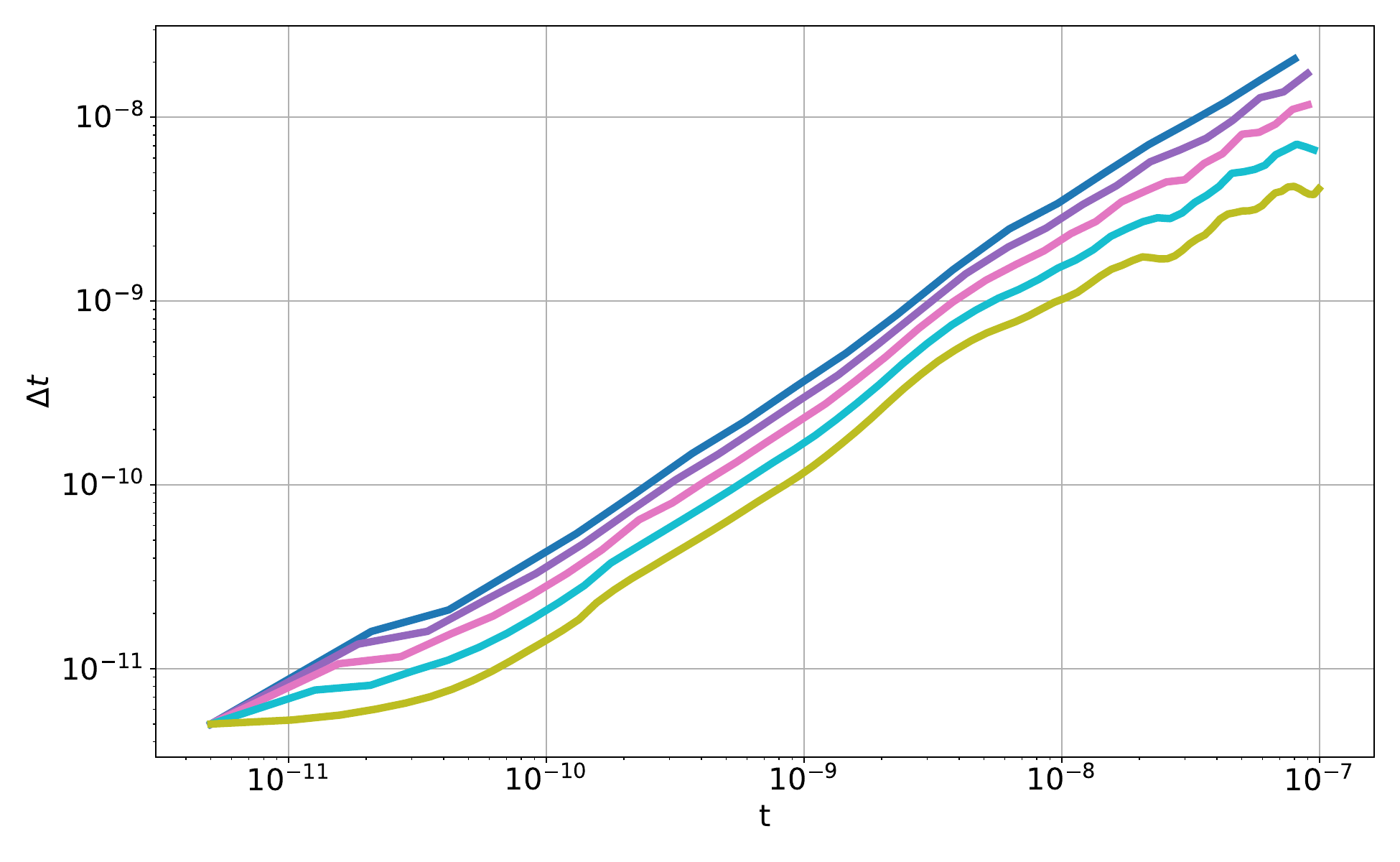}
            \caption*{(d) Semi-implicit, $T$-adaptivity}
        \end{minipage}
    \end{minipage}
    \begin{minipage}{0.18\textwidth} 
        \vspace{-2.8in}
        \centering
        \includegraphics[width=\textwidth]{figures_t_vs_dt_legend_clipped.png}
    \end{minipage}

    \caption{Plot of $\Delta t$ vs. $t$ for IMEX-NPRK[42]b using IMEX or semi-implicit integration, and $E_r$ or $T$ for error estimation.}
    \label{fig:IMEX_dt_v_t_plots_tophat}
\end{figure}

\Cref{fig:IMEX_dt_v_t_plots_tophat} then compares the adaptive timestep chosen for a fixed scheme IMEX-NPRK[42]b (similar results hold for other schemes as well) across IMEX or semi-implicit integration and the use of temperature or radiation energy for local error estimation. Here we see that, interestingly, when using temperature for error estimation rather than radiation energy, the adaptive process does not require a decrease in timestep when the radiation energy and temperature first turn the corner of the pipe. In addition, we see that the timesteps chosen for IMEX and semi-implicit integration are very similar through most of the time domain, only differing significantly at times $>10^{-8}$s. Here, the physics become truly diffusive and slow moving, and the coupled semi-implicit integration affords larger timesteps than IMEX without significant loss in accuracy due to the improved joint stability of semi-implicit integration. On a side note, this further supports the concept of decoupled IMEX integration, as error-based adaptive time stepping indicates that coupled semi-implicit integration should still take time steps for which IMEX integration is also stable.

\begin{figure}[!htbp]
    \centering
    \begin{minipage}[b]{0.49\textwidth}
        \centering
        \large{\bf Temperature}
        \includegraphics[trim={4cm 1.75cm 4cm 2.75cm},clip,width=0.78\textwidth]{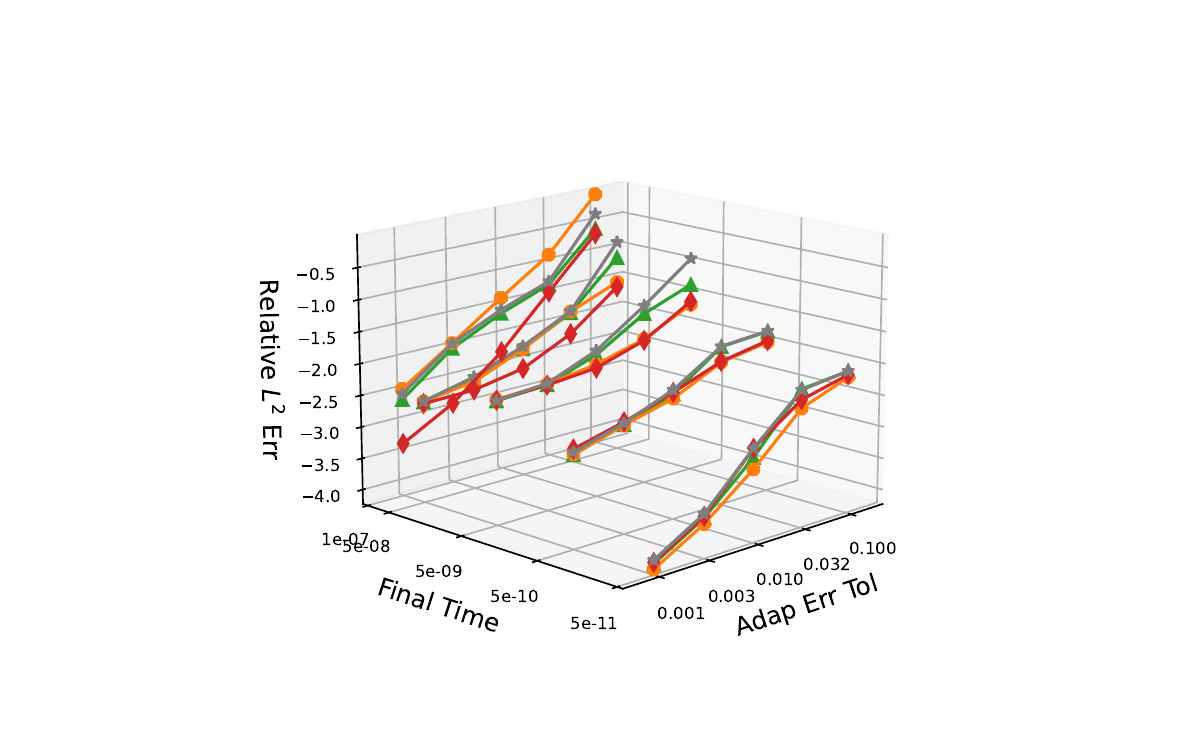} 
        \caption*{(a) IMEX, $E_r$-adaptivity}
        \hfill
        \begin{minipage}{0.3\textwidth} 
            \vspace{-4.3in}
            \centering
            \includegraphics[width=\textwidth]{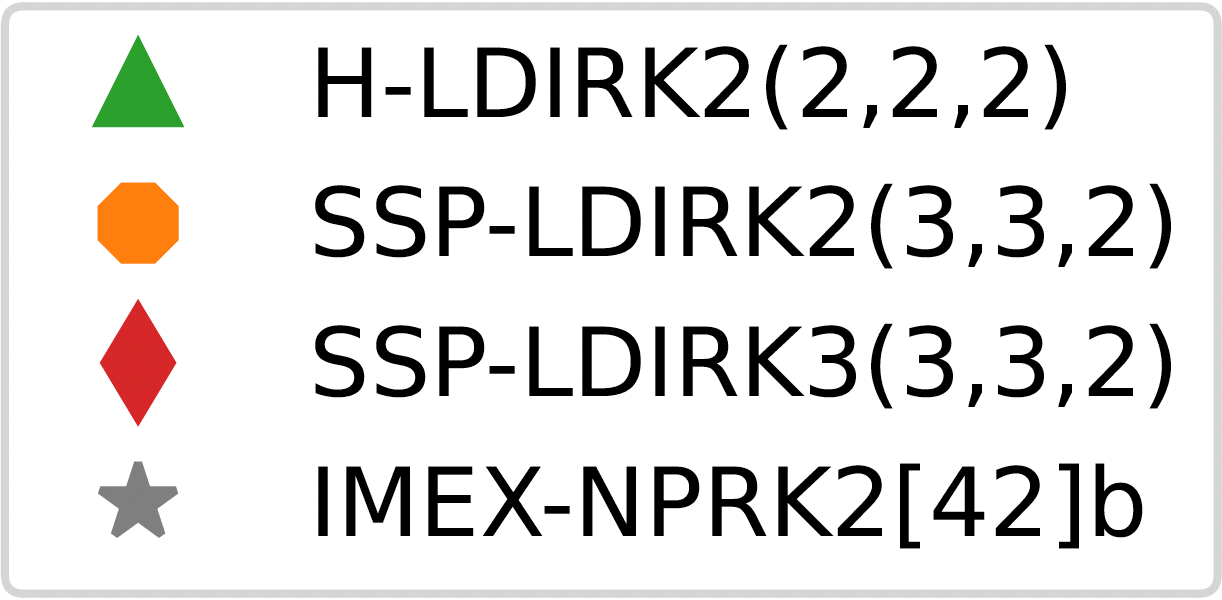}
        \end{minipage}
    \end{minipage}
    \begin{minipage}[b]{0.49\textwidth}
        \centering
        \large{\bf Radiation energy}
        \includegraphics[trim={4cm 1.75cm 4cm 2.75cm},clip,width=0.78\textwidth]{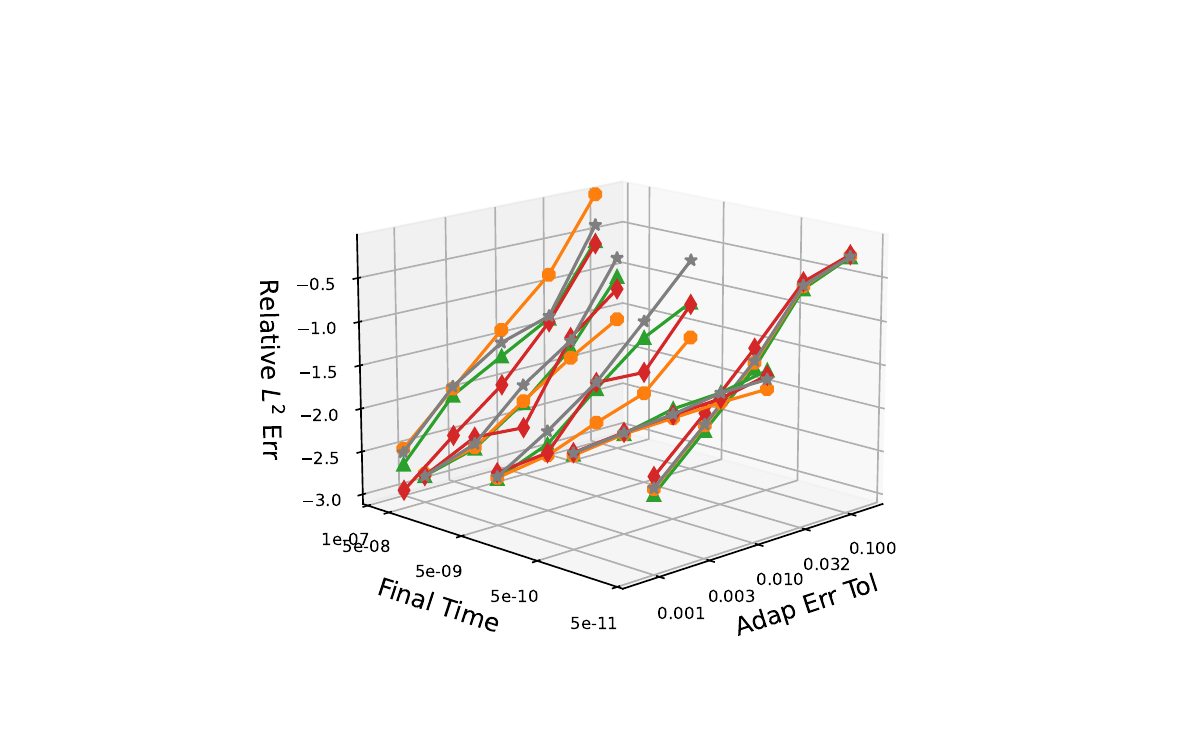} 
        \caption*{(b) IMEX, $E_r$-adaptivity}
    \end{minipage}

    \begin{minipage}[b]{0.49\textwidth}
        \centering
        \includegraphics[trim={4cm 1.75cm 4cm 2.75cm},clip,width=0.78\textwidth]{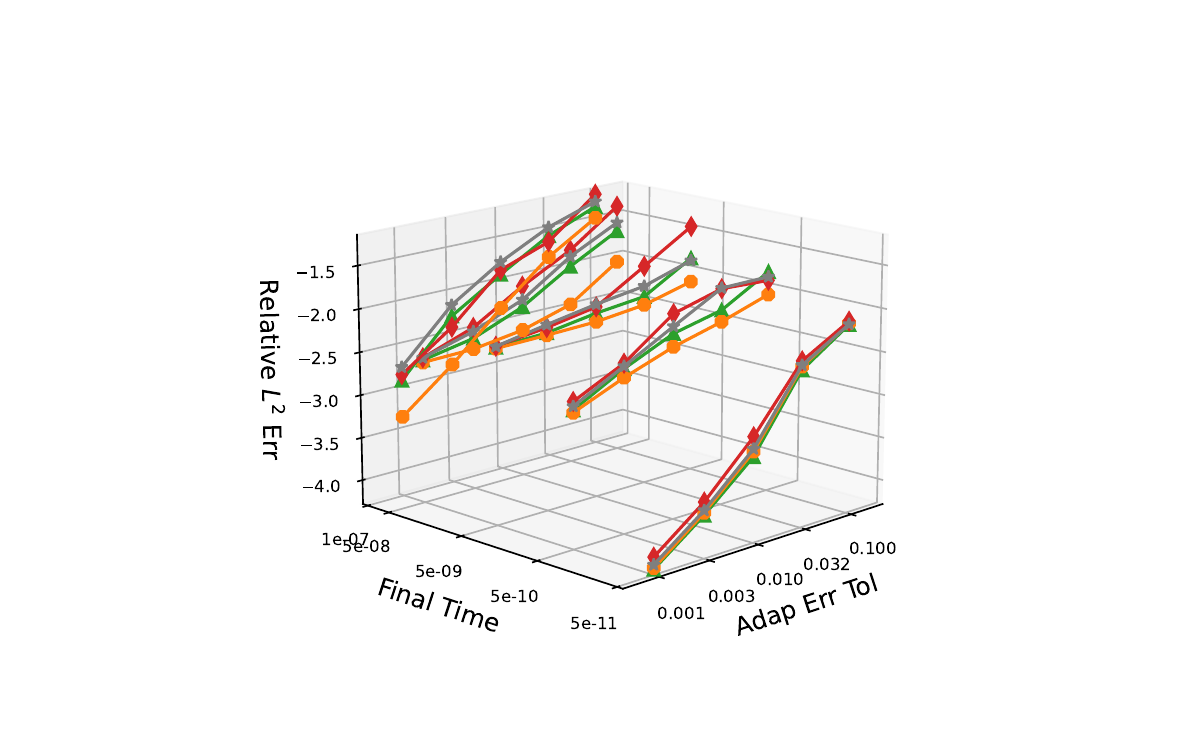} 
        \caption*{(c) Semi-implicit, $E_r$-adaptivity}
    \end{minipage}
    \begin{minipage}[b]{0.49\textwidth}
        \centering
        \includegraphics[trim={4cm 1.75cm 4cm 2.75cm},clip,width=0.78\textwidth]{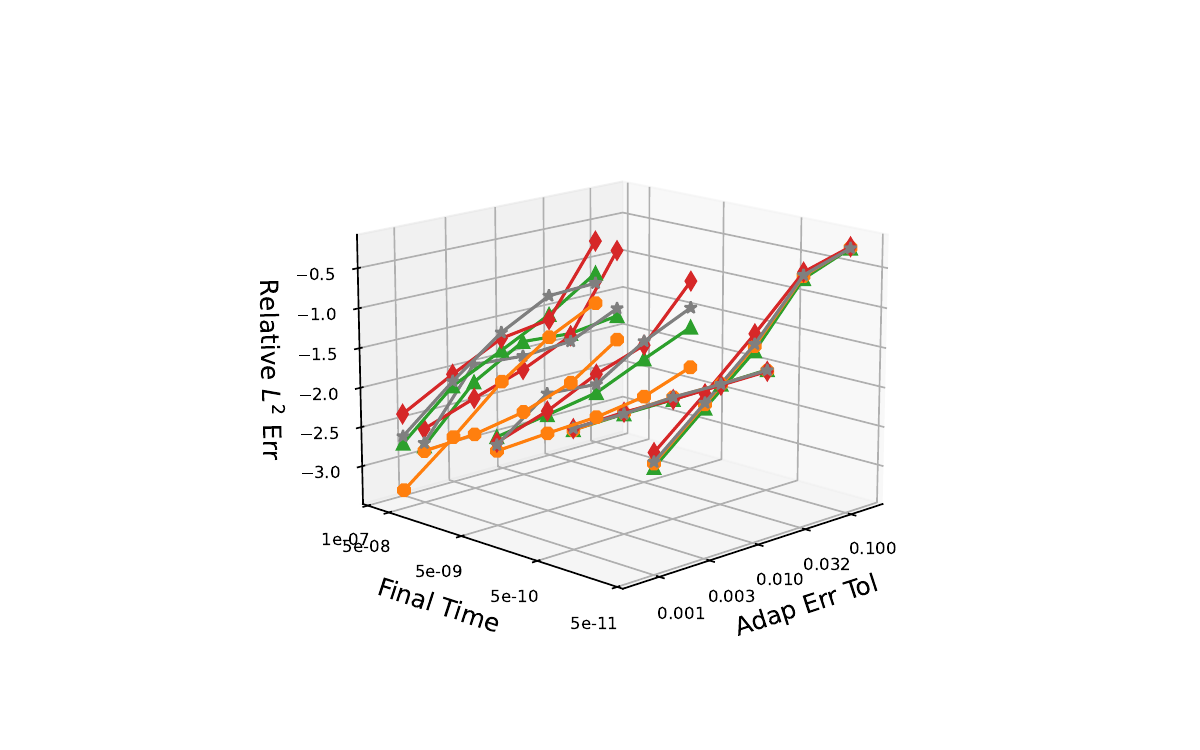} 
        \caption*{(d) Semi-implicit, $E_r$-adaptivity}
    \end{minipage}

    \begin{minipage}[b]{0.49\textwidth}
        \centering
        \includegraphics[trim={4cm 1.75cm 4cm 2.75cm},clip,width=0.78\textwidth]{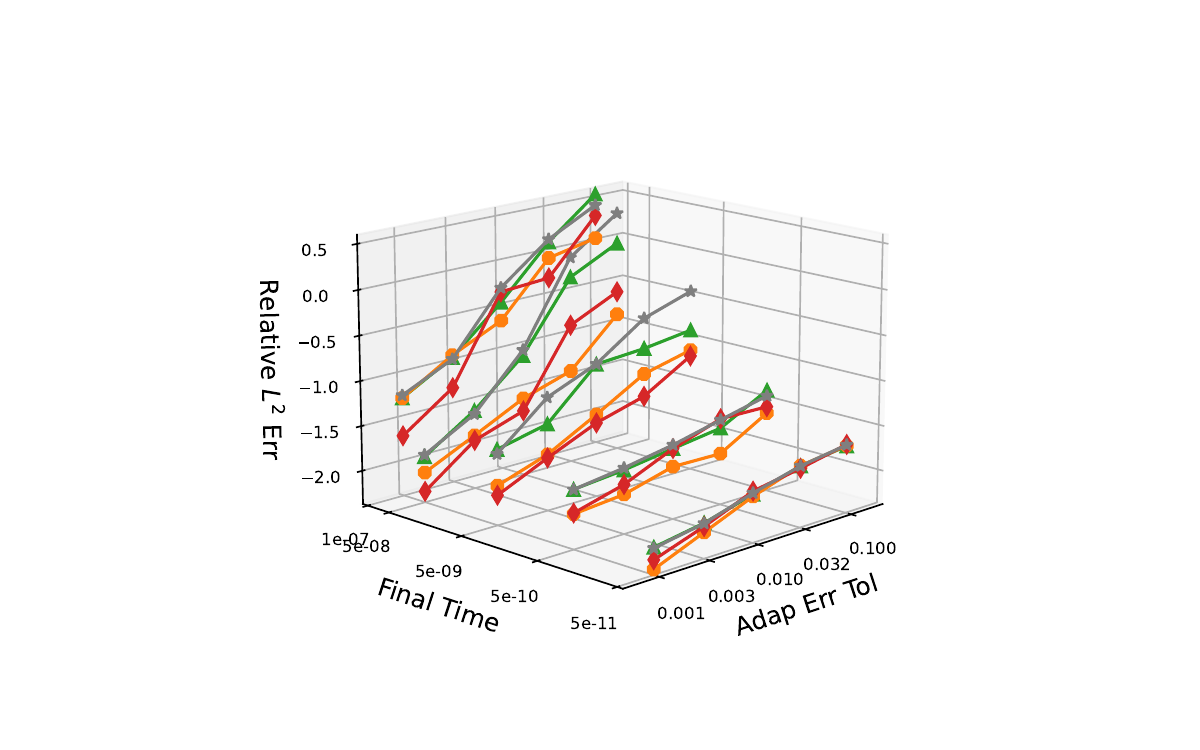} 
        \caption*{(e) IMEX, $T$-adaptivity}
    \end{minipage}
    \begin{minipage}[b]{0.49\textwidth}
        \centering
        \includegraphics[trim={4cm 1.75cm 4cm 2.75cm},clip,width=0.78\textwidth]{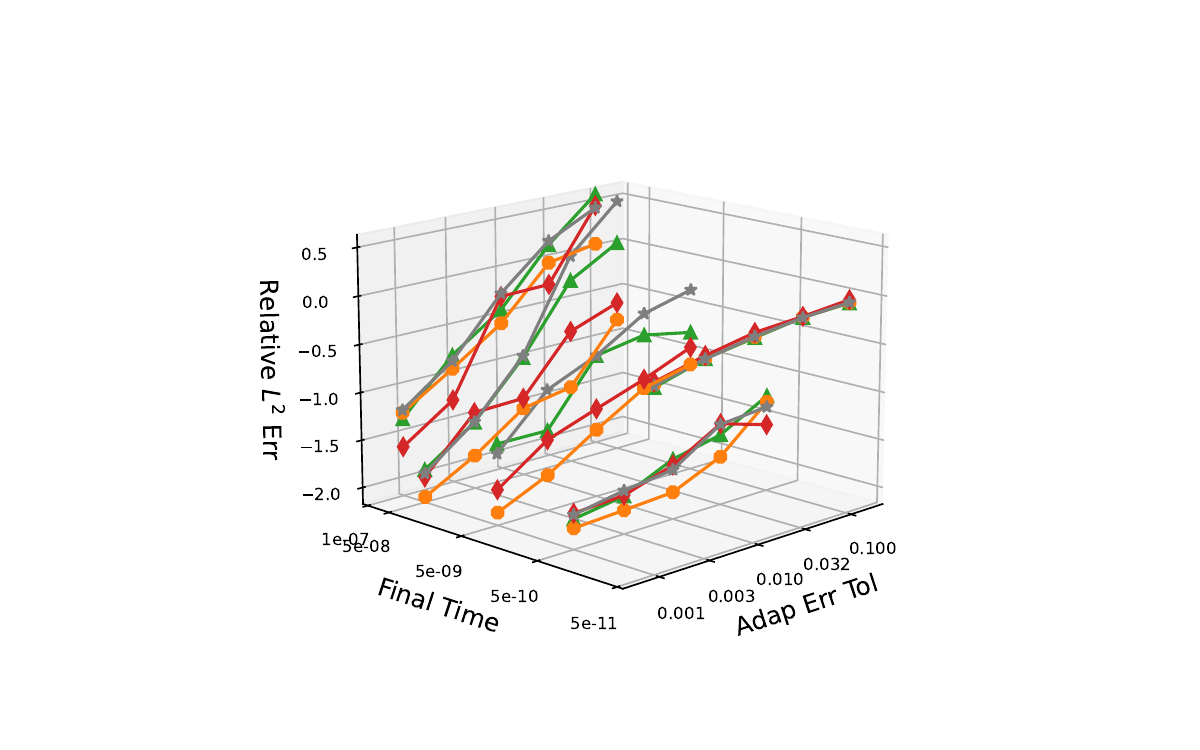} 
        \caption*{(f) IMEX, $T$-adaptivity}
    \end{minipage}

    \begin{minipage}[b]{0.49\textwidth}
        \centering
        \includegraphics[trim={4cm 1.75cm 4cm 2.75cm},clip,width=0.78\textwidth]{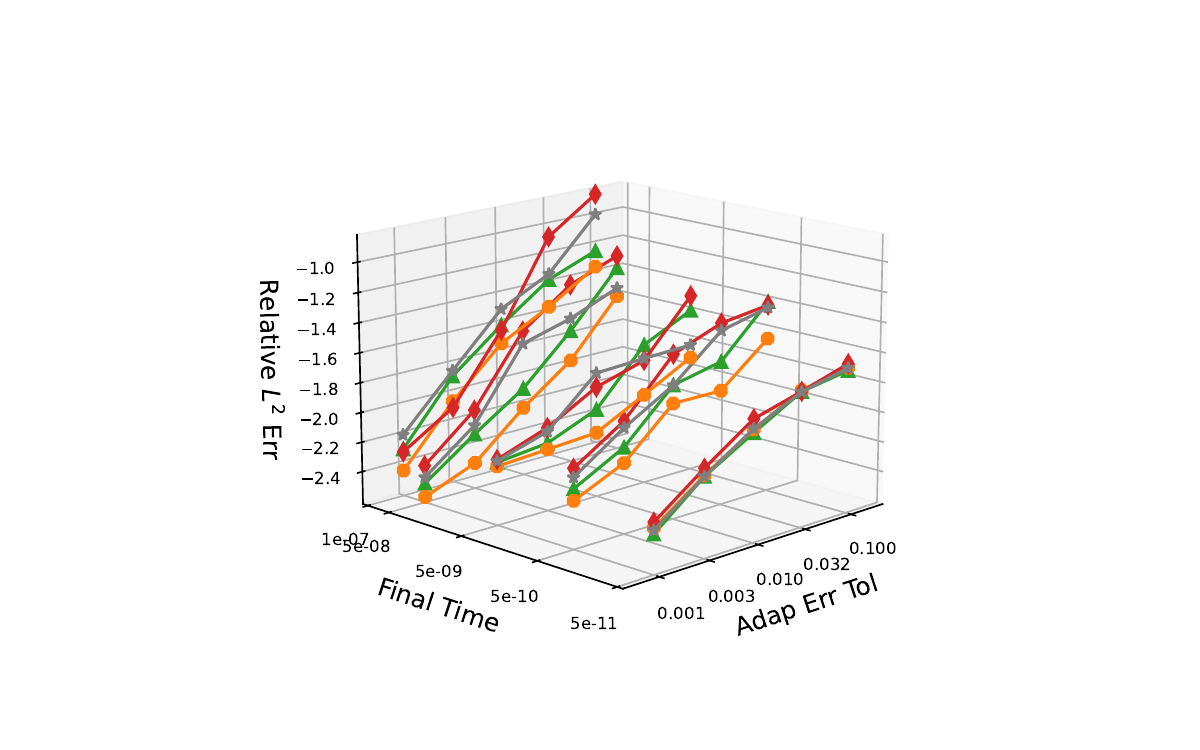} 
        \caption*{(g) Semi-implicit, $T$-adaptivity}
    \end{minipage}
    \begin{minipage}[b]{0.49\textwidth}
        \centering
        \includegraphics[trim={4cm 1.75cm 4cm 2.75cm},clip,width=0.78\textwidth]{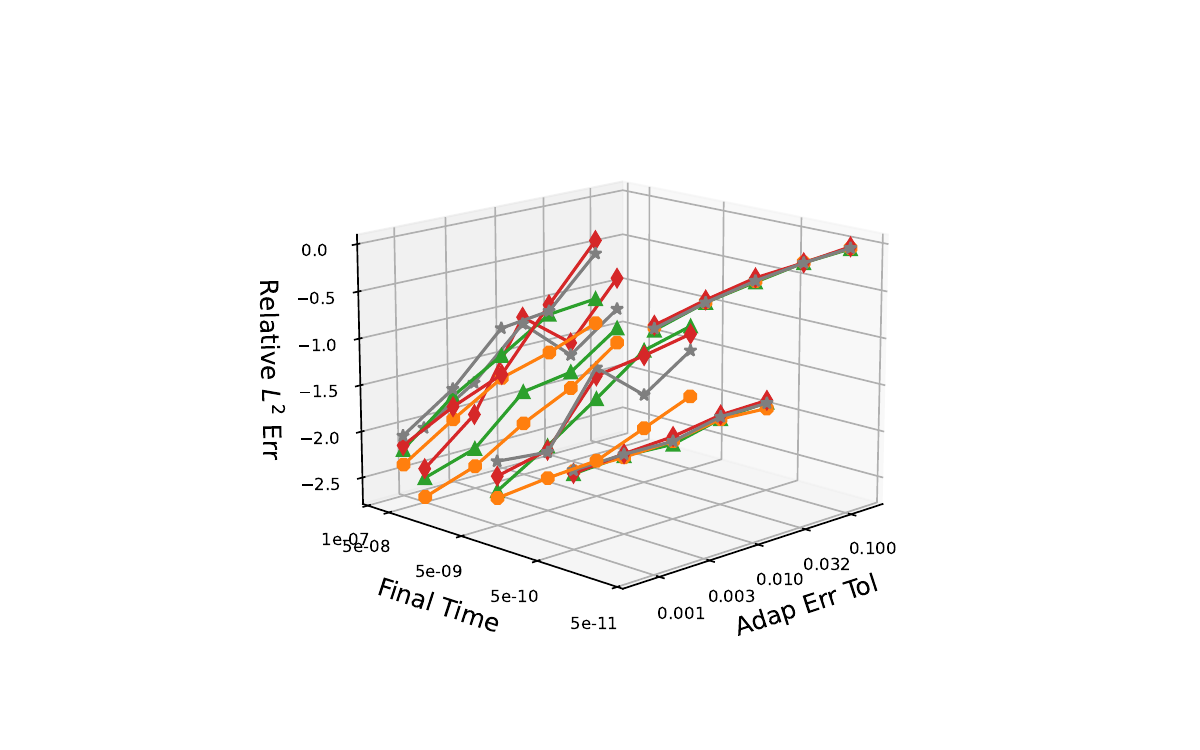} 
        \caption*{(h) Semi-implicit, $T$-adaptivity}
    \end{minipage}
    \caption{Relative $L^2$-error in log$_{10}$ scale for Tophat problem using IMEX or semi-implicit integration and $E_r$ or $T$ as the measure of error for adaptive time stepping. Left column denotes error in temperature and right error in radiation energy.}
\label{fig:convergence_tophat}
\end{figure}

Now we look at error as a function of adaptive error tolerance and final time at which the error is measured in \Cref{fig:convergence_tophat}. We consider four final times logarithmically spaced to a final time of $10^{-7}$s, at which point the dynamics are largely diffusive in nature. Because the error estimate is local in nature (in time), we demonstrate on different final times to confirm that the vastly varying time steps proposed in \Cref{fig:IMEX_dt_v_t_plots_tophat} are appropriate to accurately capture dynamics during different parts of the problem evolution, in particular transitioning from a streaming regime to a thick diffusive regime.

First, we point out that in the log-log scaling of error vs. adaptive error tolerance in \Cref{fig:convergence_tophat}, nearly all integration schemes, variables, error indicators, and final times have a clear linear relationship. That is, reducing the adaptive error tolerance by $\eta$ will reduce the integrated $L^2$-error $\sim\eta^p$ for constant $p>0$. Although $p$ varies across final time, integration scheme, variable, and error indicator, when these parameters are fixed it appears to stay relatively constant. This demonstrates that the error estimator and corresponding adaptive integration algorithm provide qualitatively accurate error estimates and refinement criteria, naturally increasing solution accuracy by reducing adaptive error tolerance.

The least successful adaptivity corresponds to early time adaptivity using $T$ to estimate error (\Cref{fig:convergence_tophat}(e-h)), where the constant $p$ is quite small, i.e. reducing the adaptive error tolerance on $T$ causes only small reductions in the observed $L^2$-error. As observed in \Cref{fig:IMEX_dt_v_t_plots_tophat}, using only $T$ to estimate error maintains a fairly large time-step early on, and it appears that $T$ alone as an error estimate has a limit in the streaming regime.  

Broadly we see that adaptivity based on $E_r$ results in smaller error than equivalent adaptivity based on $T$. This is as expected, given that $E_r$ adaptivity leads to smaller adaptive time step selection (\Cref{fig:IMEX_dt_v_t_plots_tophat}). Similarly, semi-implicit integration leads to more accurate results than IMEX, particularly in late time when the problem is diffusive in nature. This is also expected, since the adaptive timesteps proposed are quite similar between IMEX and semi-implicit (\Cref{fig:IMEX_dt_v_t_plots_tophat}), but semi-implicit integration better resolves the implicit coupling between high-order transport phenomena and low-order diffusion. Last, we point out that the error obtained is in most cases quite similar across all four schemes. This supports the use of the two-stage schemes, H-LDIRK2(2,2,2) and IMEX-NPRK2[42]b, in practice due to their lower computational costs and storage requirements.

\subsection{Larsen}\label{sec:results:larsen}

\begin{figure}[!ht]
\centering

\begin{subfigure}{0.3\textwidth}
  \centering
  \includegraphics[width=\textwidth]{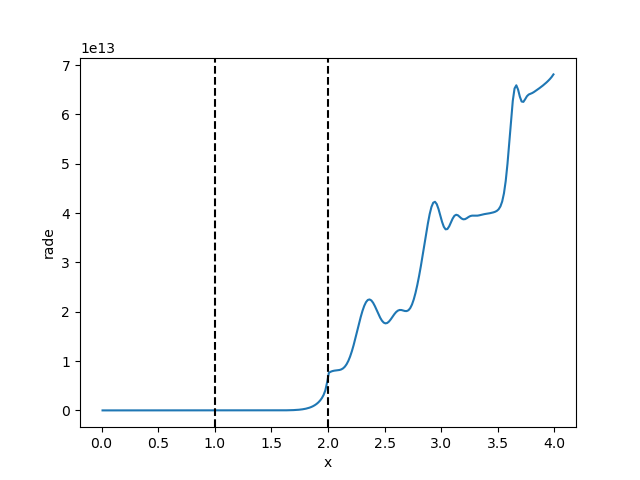}
  \caption{}
\end{subfigure}
\hfill
\begin{subfigure}{0.3\textwidth}
  \centering
  \includegraphics[width=\textwidth]{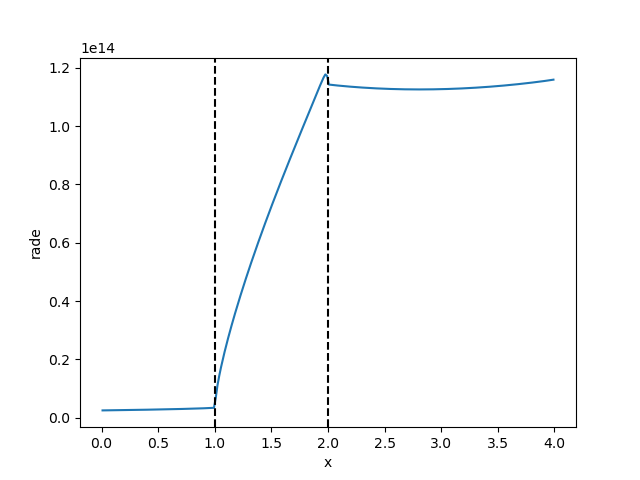}
  \caption{}
\end{subfigure}
\hfill
\begin{subfigure}{0.3\textwidth}
  \centering
  \includegraphics[width=\textwidth]{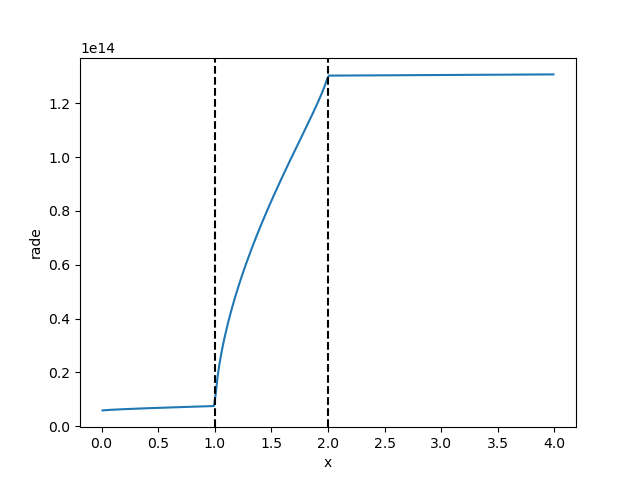}
  \caption{}
\end{subfigure}

\vspace{1em}

\begin{subfigure}{0.3\textwidth}
  \centering
  \includegraphics[width=\textwidth]{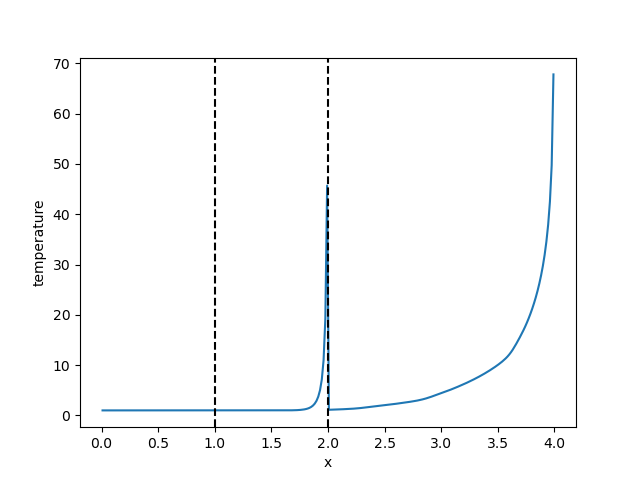}
  \caption{}
\end{subfigure}
\hfill
\begin{subfigure}{0.3\textwidth}
  \centering
  \includegraphics[width=\textwidth]{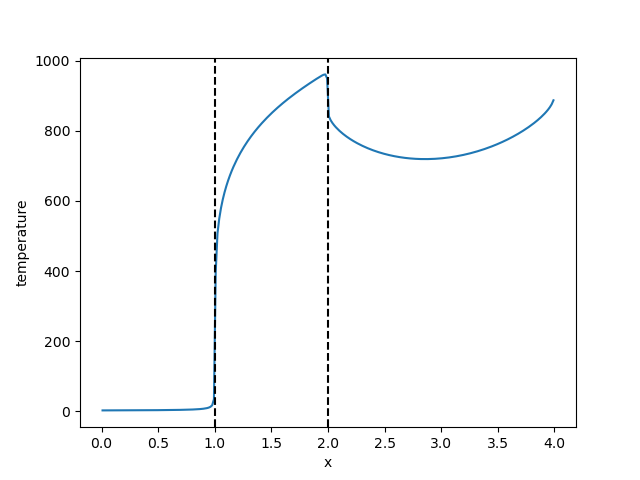}
  \caption{}
\end{subfigure}
\hfill
\begin{subfigure}{0.3\textwidth}
  \centering
  \includegraphics[width=\textwidth]{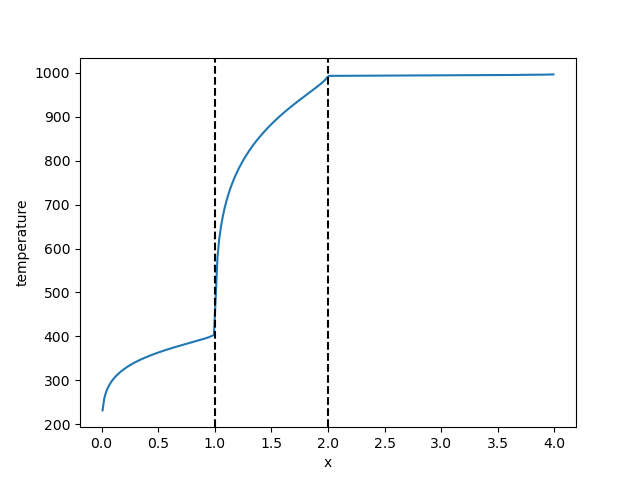}
  \caption{}
\end{subfigure}

\caption{{\bf{Top:}} Radiation energy. {\bf{Bottom:}} Material temperature. Shown are solution snapshots at $t=7.64968\times 10^{-11}, 2.88063\times 10^{-9}$, and $10^{-7}$ corresponding to the two drops in $\Delta t$ seen in Figure \ref{fig:IMEX_dt_v_t_plots_larsen} and the final time, respectively}
\label{fig:larsen_temp_and_rade_snapshots}
\end{figure}

The Larsen problem is specified as a 1d multifrequency problem from \cite[Sec. VII]{Larsen.1988} with constant heat capacity $C_V = 5.109\cdot 10^{11}$ \si{\erg\per\eV\per\cm\cubed} and opacity specified by
\begin{equation*}
    \sigma(\nu,T) \coloneqq \frac{\gamma}{(h\nu)^3}\left(1-e^{-h\nu/T}\right),
\end{equation*}
where $h$ is Planck's constant and $\gamma$ a specified constant. For 1d domain $[0,4]$cm, we have 
\begin{equation*}
    \gamma(x) = \begin{cases}
        10^9 \textnormal{ eV$^3$/cm} & 0<x<1 \\
        10^{12} \textnormal{ eV$^3$/cm} & 1<x<2 \\
        10^9 \textnormal{ eV$^3$/cm} & 2<x<4
    \end{cases}.
\end{equation*}
and utilize inflow boundary conditions on the right according to a Planckian distribution evaluated at a temperature of $T_b = $ \SI{1000}{\eV}, that is $I_\text{inflow} = \frac{ac T_b^4}{4\pi}$. The initial temperature and radiation fields are set to be in equilibrium with each other at a spatially uniform temperature of \SI{1}{\eV}. Here we present analogous results for the Larsen problem as demonstrated for the Tophat problem in \Cref{sec:results:tophat}. We run this problem on a uniform grid of $256$ cells with 8th-order S$_N$ quadrature, 50 energy groups uniformly split in log-scale between $10^{-2}$eV and $10^6$eV., and a HOLO discretization as in \cite{Park.2020,imex-trt}. A sequence of solutions for radiation energy and temperature are shown in \Cref{fig:larsen_temp_and_rade_snapshots}.

\begin{figure}[!htbp]
    \centering
    \begin{minipage}{0.75\textwidth} 
        \centering
        \begin{minipage}[b]{0.48\textwidth} 
            \centering
            \includegraphics[width=\textwidth]{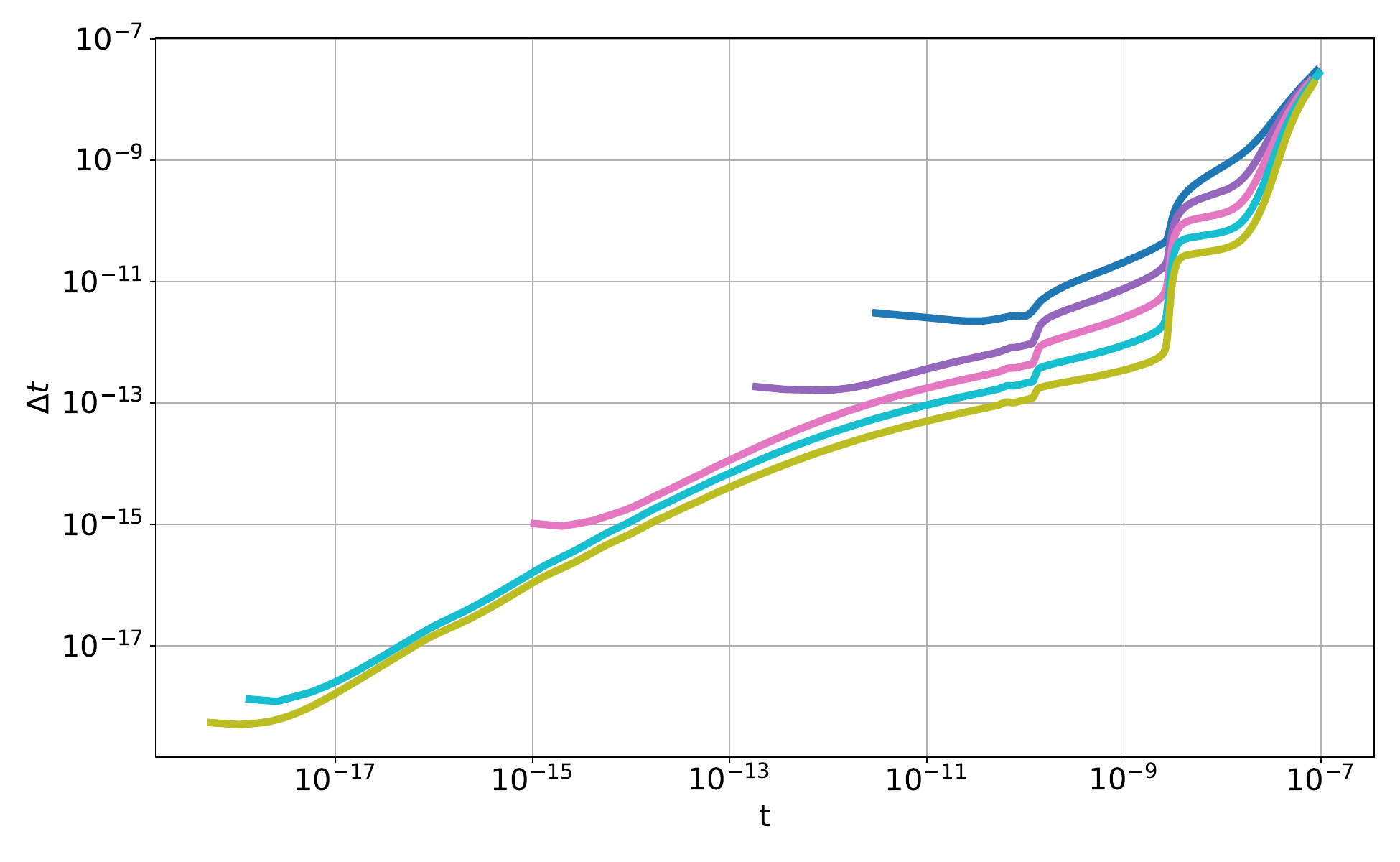}
            \caption*{(a) IMEX, $E_r$-adaptivity}
        \end{minipage}
        \hfill
        \begin{minipage}[b]{0.48\textwidth} 
            \centering
            \includegraphics[width=\textwidth]{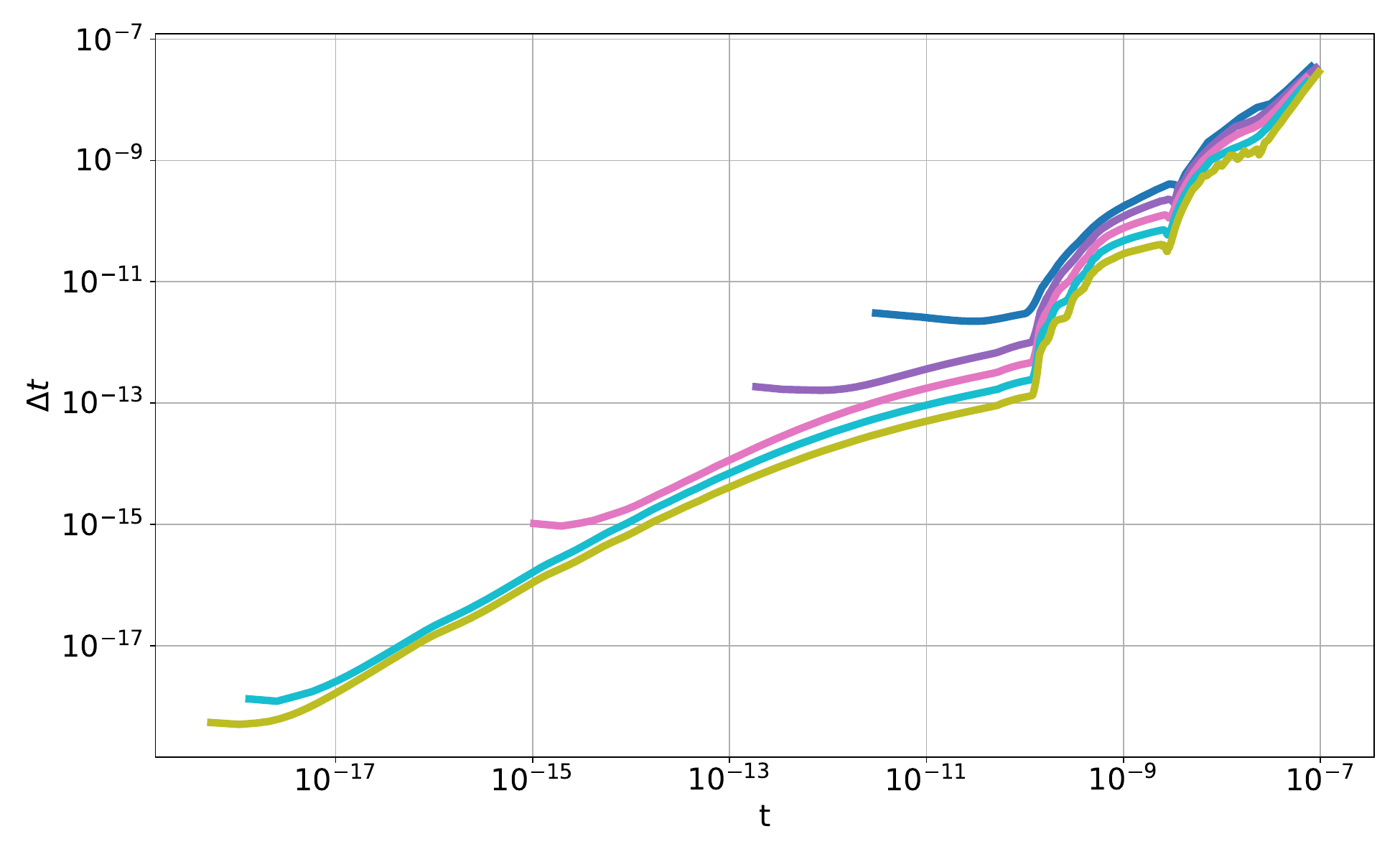}
            \caption*{(b) Semi-implicit, $E_r$-adaptivity}
        \end{minipage}

        \vspace{5pt} 

        \begin{minipage}[b]{0.48\textwidth} 
            \centering
            \includegraphics[width=\textwidth]{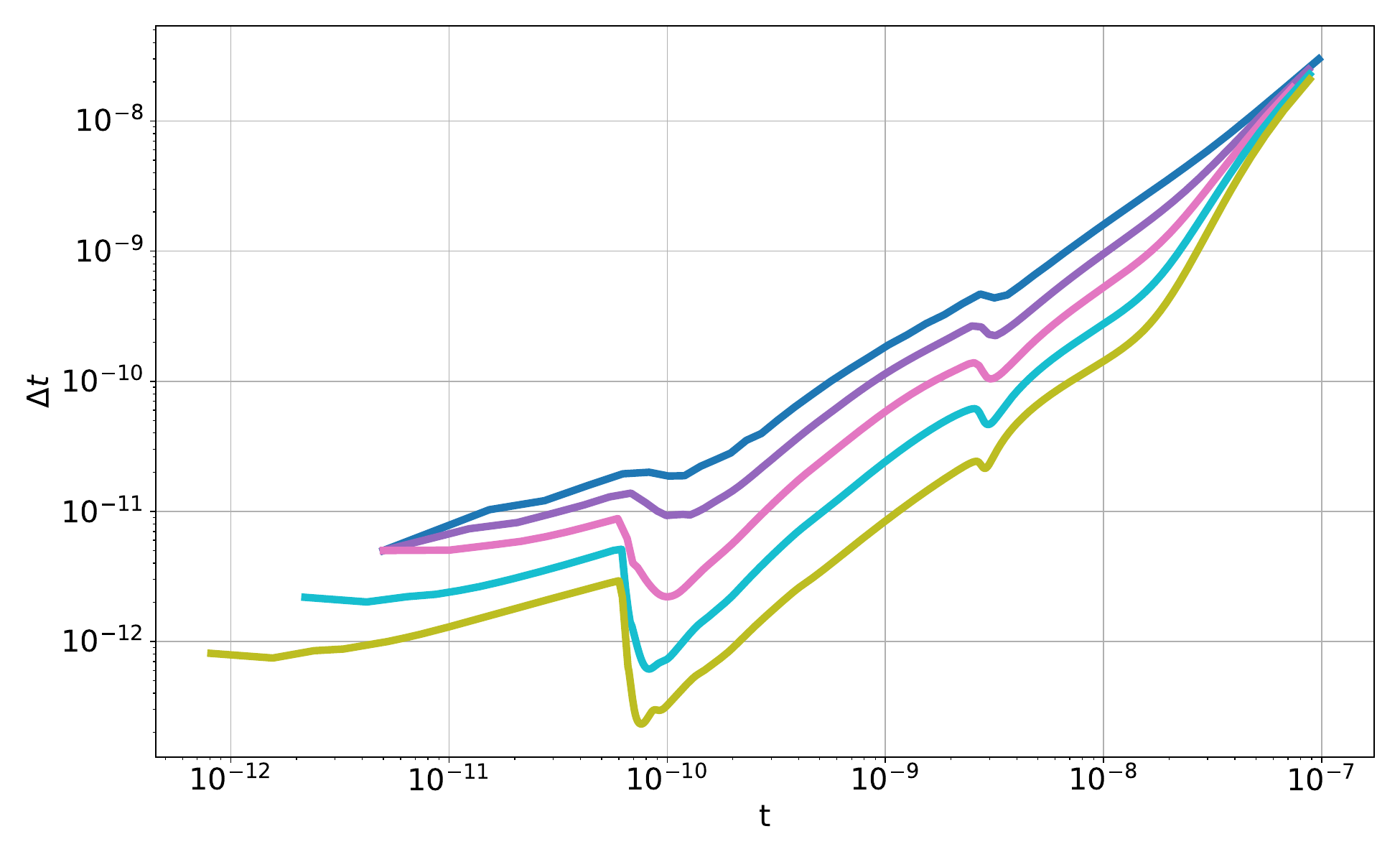}
            \caption*{(c) IMEX, $T$-adaptivity}
        \end{minipage}
        \hfill
        \begin{minipage}[b]{0.48\textwidth} 
            \centering
            \includegraphics[width=\textwidth]{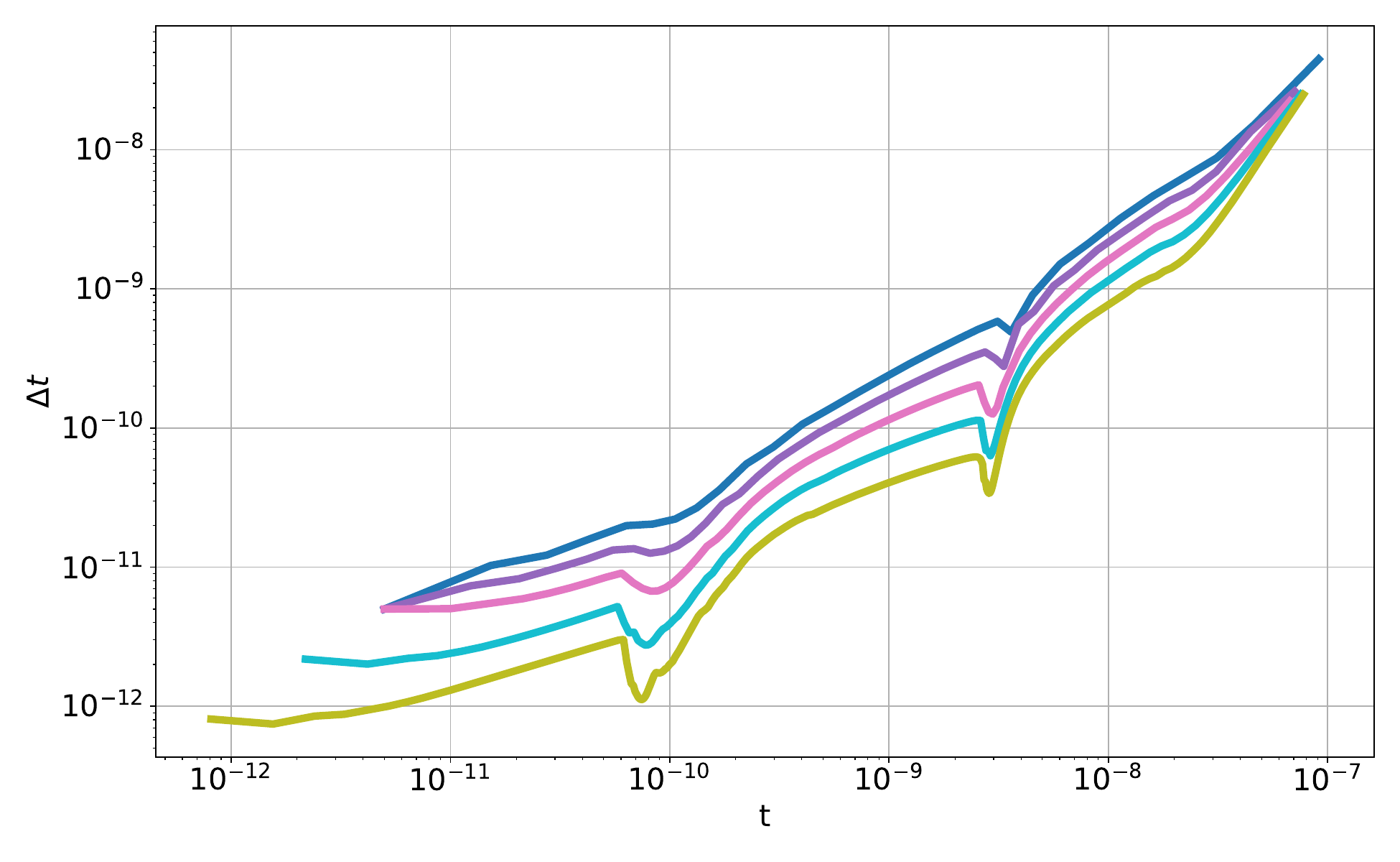}
            \caption*{(d) Semi-implicit, $T$-adaptivity}
        \end{minipage}
    \end{minipage}
    \begin{minipage}{0.18\textwidth} 
        \vspace{-2.8in}
        \centering
        \includegraphics[width=\textwidth]{figures_t_vs_dt_legend_clipped.png}
    \end{minipage}

    \caption{Plot of $\Delta t$ vs. $t$ for IMEX-NPRK[42]b applied to the Larsen problem using IMEX or semi-implicit integration, and $E_r$ or $T$ for error estimation.}
    \label{fig:IMEX_dt_v_t_plots_larsen}
\end{figure}

We start by again presenting adaptive time step as a function of physical time for four different local error tolerances, and for IMEX or semi-implicit integration and using temperature or radiation energy for local error estimation in \Cref{fig:IMEX_dt_v_t_plots_larsen}. Here, the time steps arising from the four different IMEX schemes and embedded pairs are nearly identical (in contrast to the tophat problem, where SSP-LDIRK2(3,3,2) had an instability of sorts), so we only present results for IMEX-NPRK[42]b. As for the tophat problem, we see that error estimation based on $E_r$ drives the timestep significantly smaller than error estimation based on $T$. Looking at the temperature adaptivity in \Cref{fig:IMEX_dt_v_t_plots_larsen}, we see two sudden drops in $\Delta t$, each also corresponding to a sharp increase in $\Delta t$ using $E_r$-adaptivity. Thus in both cases, these two times are marked as having a significant change in dynamic timescale. As expected, these times correspond to the radiation energy heating up a new material. \Cref{fig:larsen_temp_and_rade_snapshots}(a,d) show the solution shortly after the large drop in $\Delta t$ in \Cref{fig:IMEX_dt_v_t_plots_larsen}, where the second, optically thicker material begins to heat up. Similarly, the smaller drop in $\Delta t$ seen in \Cref{fig:IMEX_dt_v_t_plots_larsen} corresponds to the dynamics seen in \Cref{fig:larsen_temp_and_rade_snapshots}(b,e), depicting the initial increase in temperature of the third material. 

\begin{figure}[htbp]
    \centering
    \begin{minipage}[b]{0.49\textwidth}
        \centering
        \large{\bf Temperature}
        \includegraphics[trim={4cm 1.75cm 4cm 2.75cm},clip,width=0.78\textwidth]{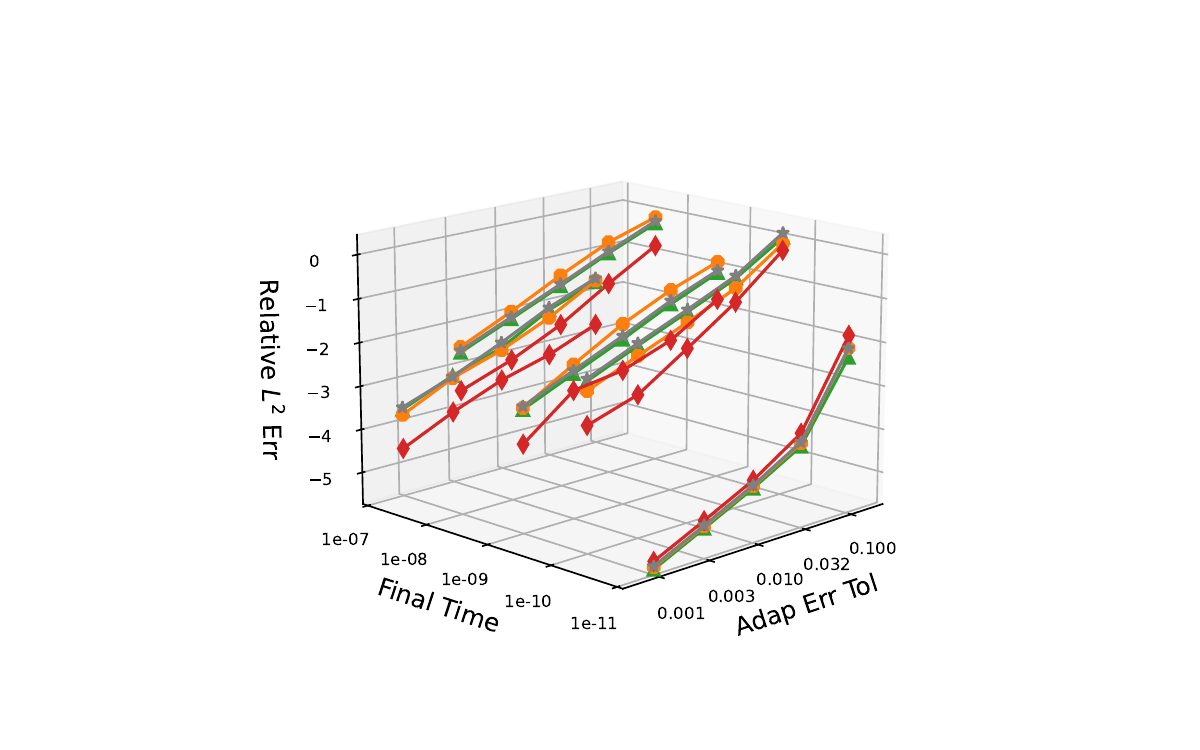} 
        \caption*{(a) IMEX, $E_r$-adaptivity}
        \hfill
    \begin{minipage}{0.3\textwidth} 
            \vspace{-4.2in}
            \centering
            \includegraphics[width=\textwidth]{figures_conv_plots_legend_clipped.png}
        \end{minipage}
    \end{minipage}
    \begin{minipage}[b]{0.49\textwidth}
        \centering
        \large{\bf Radiation energy}
        \includegraphics[trim={4cm 1.75cm 4cm 2.75cm},clip,width=0.78\textwidth]{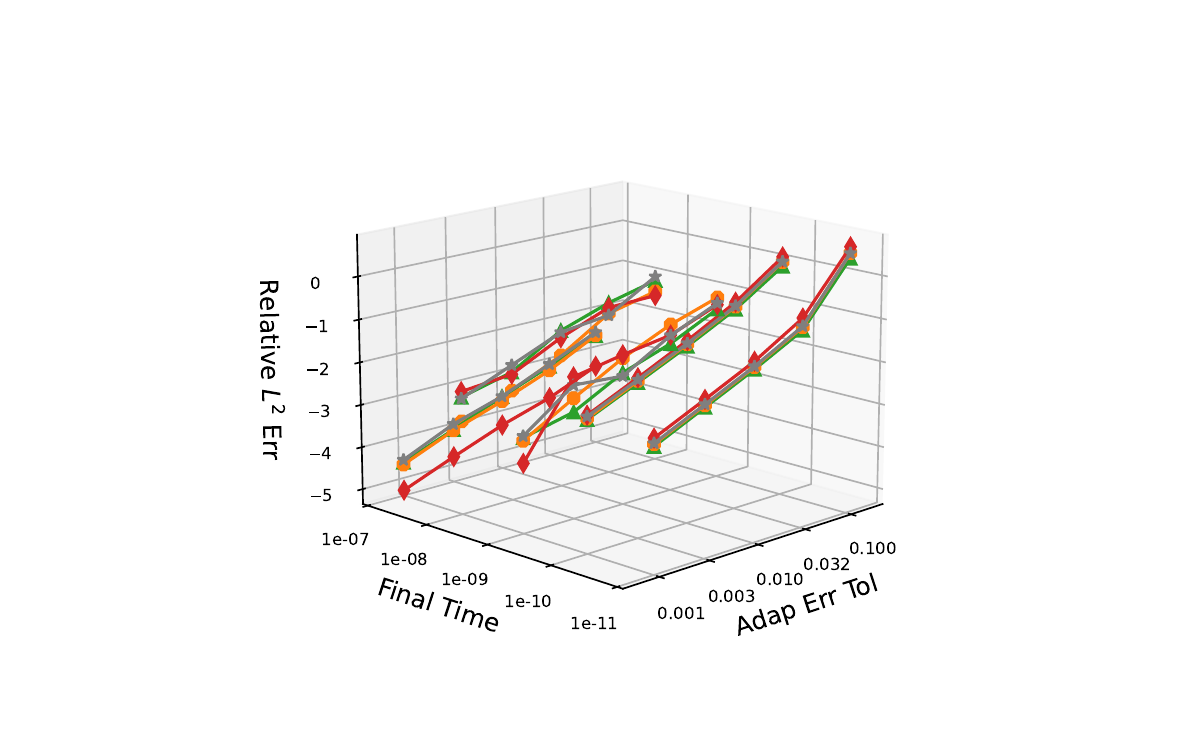} 
        \caption*{(b) IMEX, $E_r$-adaptivity}
    \end{minipage}

    \vspace{1ex}
    \begin{minipage}[b]{0.49\textwidth}
        \centering
        \includegraphics[trim={4cm 1.75cm 4cm 2.75cm},clip,width=0.78\textwidth]{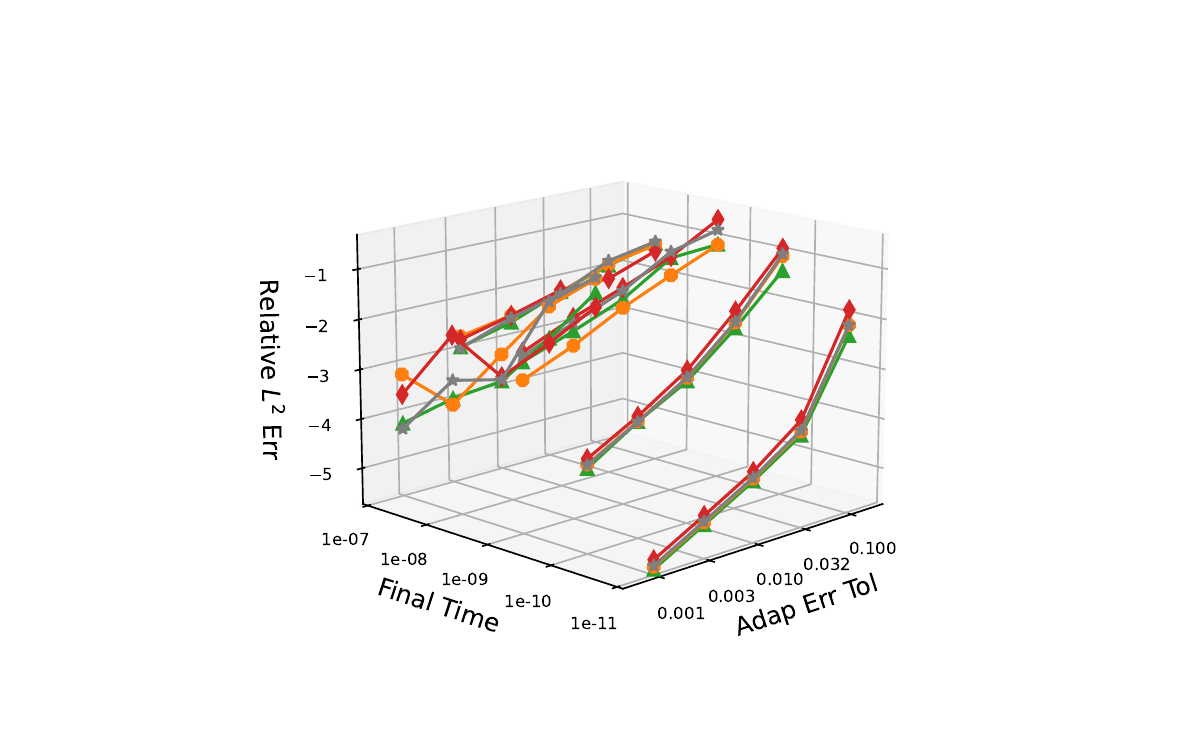} 
        \caption*{(c) Semi-implicit, $E_r$-adaptivity}
    \end{minipage}
    \begin{minipage}[b]{0.49\textwidth}
        \centering
        \includegraphics[trim={4cm 1.75cm 4cm 2.75cm},clip,width=0.78\textwidth]{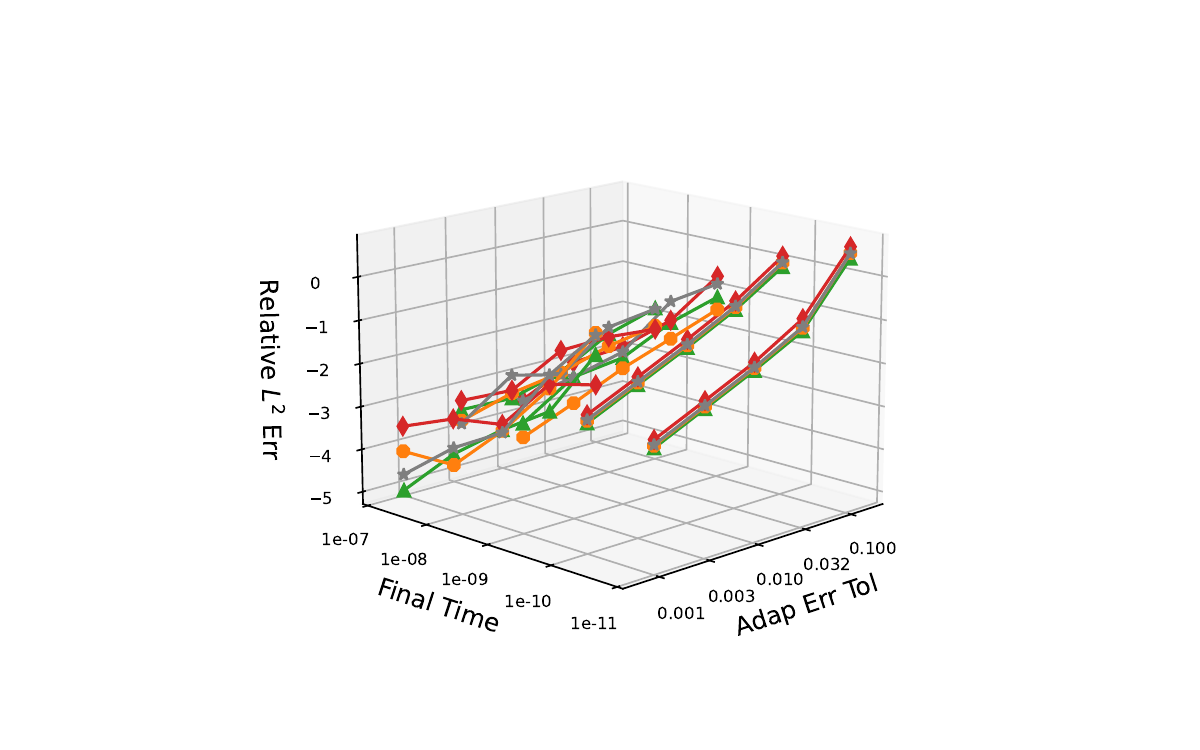} 
        \caption*{(d) Semi-implicit, $E_r$-adaptivity}
    \end{minipage}

    \begin{minipage}[b]{0.49\textwidth}
        \centering
        \includegraphics[trim={4cm 1.75cm 4cm 2.75cm},clip,width=0.78\textwidth]{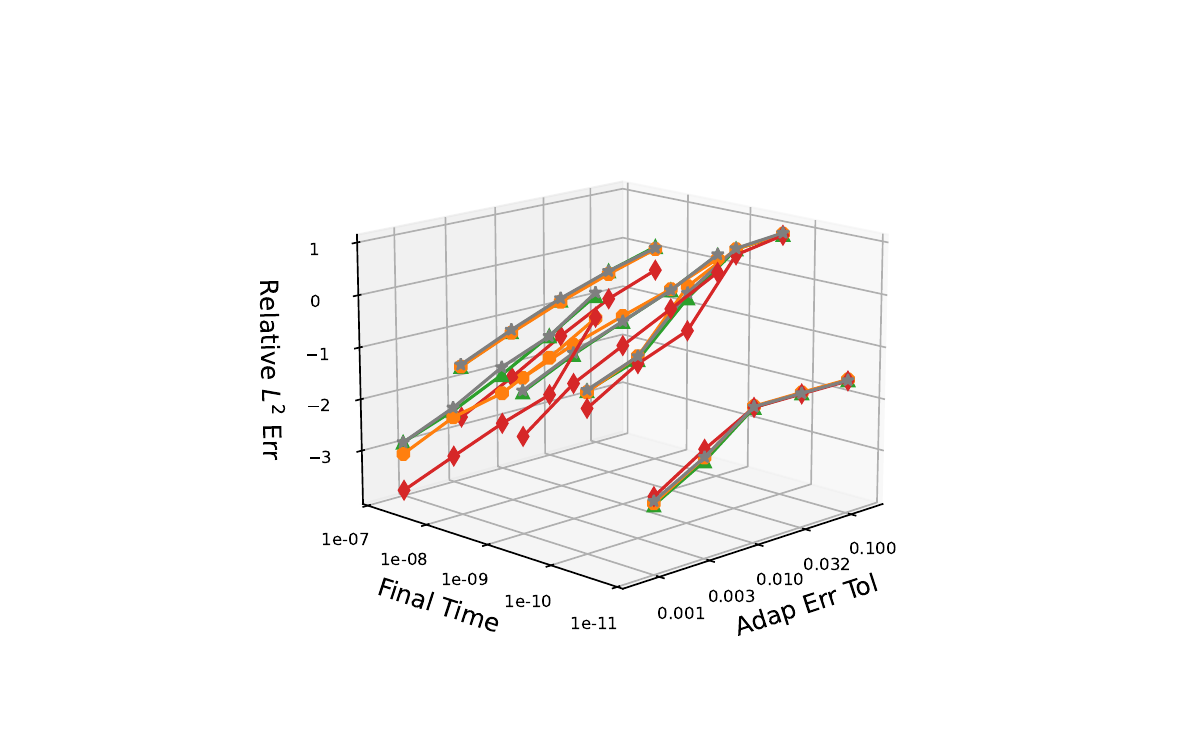} 
        \caption*{(e) IMEX, $T$-adaptivity}
    \end{minipage}
    \begin{minipage}[b]{0.49\textwidth}
        \centering
        \includegraphics[trim={4cm 1.75cm 4cm 2.75cm},clip,width=0.78\textwidth]{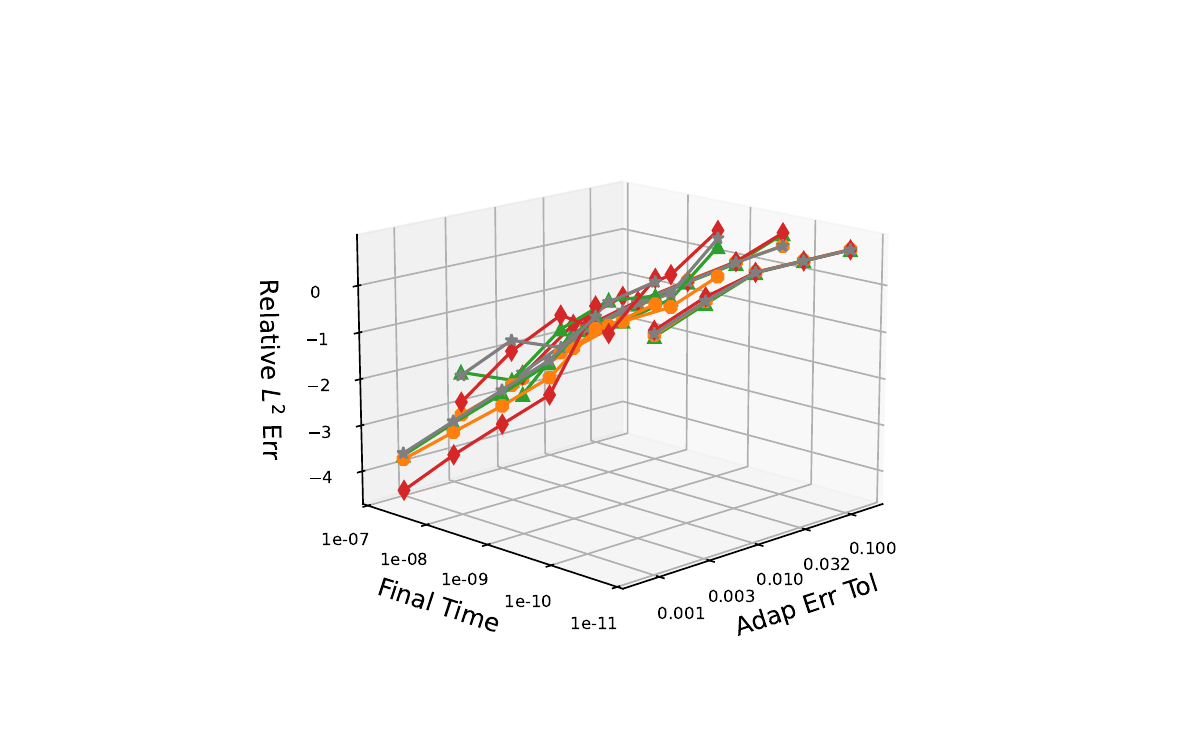} 
        \caption*{(f) IMEX, $T$-adaptivity}
    \end{minipage}

    \begin{minipage}[b]{0.49\textwidth}
        \centering
        \includegraphics[trim={4cm 1.75cm 4cm 2.75cm},clip,width=0.78\textwidth]{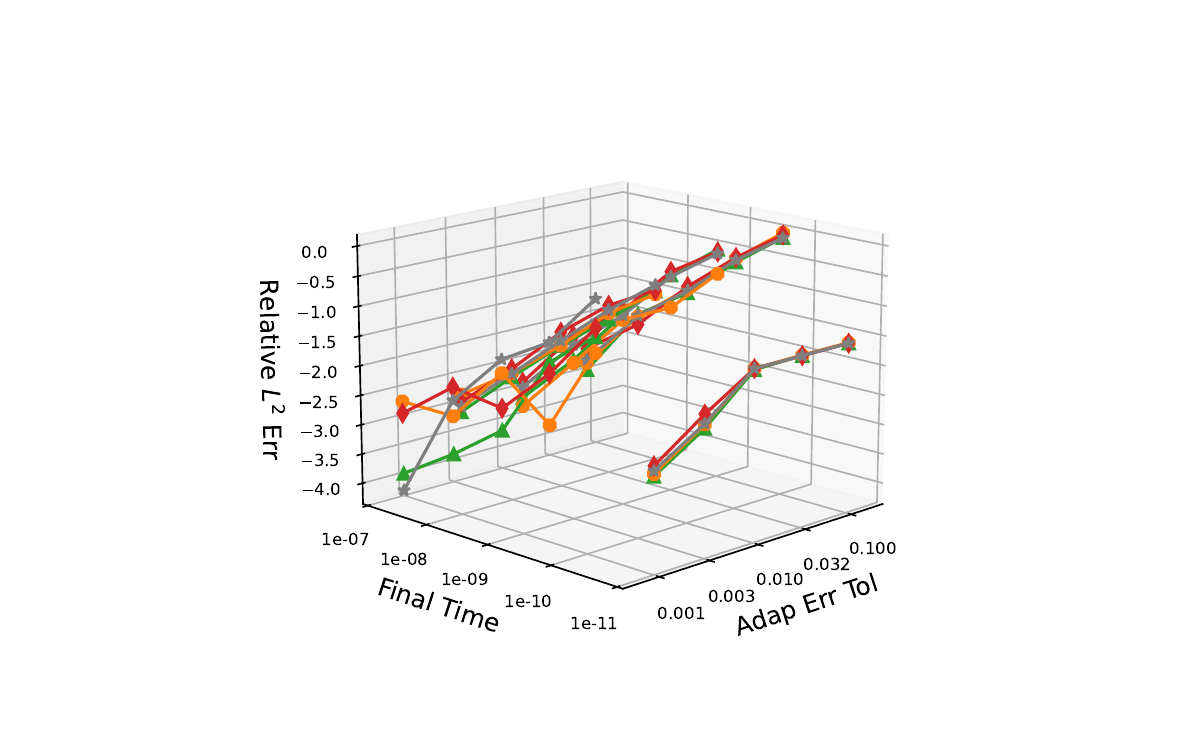} 
        \caption*{(g) Semi-implicit, $T$-adaptivity}
    \end{minipage}
    \begin{minipage}[b]{0.49\textwidth}
        \centering
        \includegraphics[trim={4cm 1.75cm 4cm 2.75cm},clip,width=0.78\textwidth]{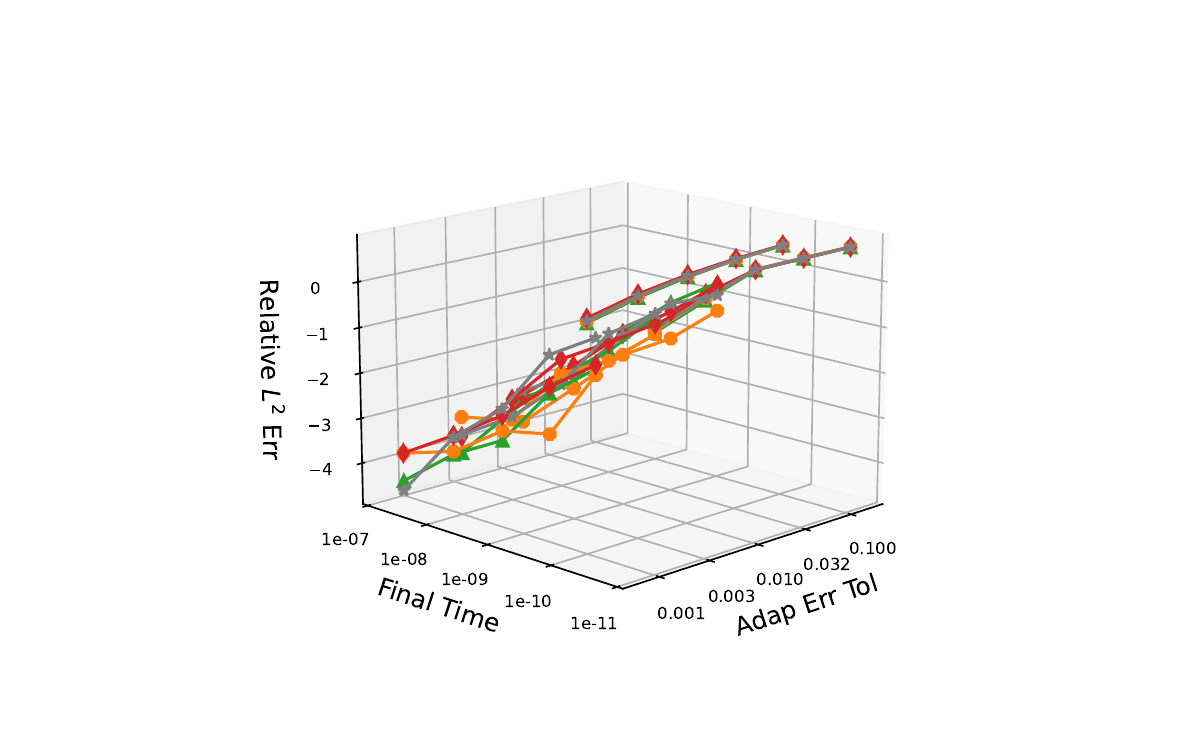} 
        \caption*{(h) Semi-implicit, $T$-adaptivity}
    \end{minipage}
    \caption{Relative $L^2$-error in log$_{10}$ scale for Larsen problem using IMEX or semi-implicit integration and $E_r$ or $T$ as the measure of error for adaptive time stepping. Left column denotes error in temperature and right error in radiation energy.}
\label{fig:convergence_larsen}
\end{figure}

Now we again look at error as a function of adaptive error tolerance and final time at which the error is measured in \Cref{fig:convergence_larsen}. 
As for tophat, the log-log scaling of error vs. adaptive error tolerance has a linear relationship in most cases, excluding the very first final time considered, where the relationship is more complex. Similarly, almost all cases yield a steady reduction in error arising from a reduced adaptive error tolerance, with the only exception being $T$-adaptivity not effectively reducing radiation energy error in very early time. Using radiation energy as the error estimator yields smaller $L^2$-error, but noting how much smaller the proposed timesteps are compared with using temperature as an error estimator (see \Cref{fig:IMEX_dt_v_t_plots_larsen}), this accuracy comes at a cost. Last, we point out that for the Larsen problem, semi-implicit integration does not offer significant benefit over IMEX integration in terms of $L^2$-error, despite very similar adaptive time step selection (\Cref{fig:IMEX_dt_v_t_plots_larsen}). This again supports the idea that the IMEX partitioning is stable and accurate on the dynamic timescale we want to capture. 

{\color{black}
Part of the novel contribution of this paper is the extension of the IMEX formulation from \cite{imex-trt} for gray TRT to multifrequency TRT problems, as well as the introduction of a high-order semi-implicit method for gray or multifrequency that is only explicit/linearized in opacity. To this end, we conclude by providing a convergence analysis for both the semi-implicit and IMEX methods with fixed $\Delta t$ for the multifrequency Larsen problem. For each of the presented integrators, we compute the solution on six uniform grids in time.  The initial timestep is chosen to be $\Delta t_0 = 3.2$e-$11$ and the timestep is halved for each refinement level so that $\Delta t_6 = 1$e-$12$.  We use the SSP-LDIRK(3,3,2) integrator to compute our reference solutions with $\Delta t_{ref}=1$e-$13$ since the SSP-LDIRK(3,3,3) scheme is consistently one of the most accurate in the previously presented convergence studies.  For this fixed-step convergence study, we set the final time to be $t_{final} = 1$e-$9$.  We provide convergence tables for $E_r$ and the temperature in Tables \ref{tab:convergence_imex_temperature}-\ref{tab:convergence_semi_rade}. In each case, we compute the relative $L^1$, $L^2$ and $L^\infty$ norms with their respective convergence rates given in parentheses.

The data shows that the IMEX method, without adaptive time stepping, suffers from order reduction at these larger time steps (compared to adaptive $\Delta t$ for early physical time), observing convergence usually in the range $\rho\in[1.2,1.5]$. This is very likely due to the more complex relationship between LO and HO opacities in the multifrequency setting compared to gray, wherein we only obtain a ``correct'' LO opacity by resolving the HO equation in frequency. Consistent with this explanation, the semi-implicit schemes that accurately resolve the frequency dependence in LO opacities (despite linearizing the temperature dependence) attain the desired convergence properties even for large time steps, as well as error 2-3 orders of magnitude smaller than the IMEX methods. What is interesting is that the adaptive time stepping provides qualitatively similar behavior in choice of $\Delta t$ for the IMEX and semi-implicit formulations (see \Cref{fig:IMEX_dt_v_t_plots_larsen}), and, moreover, the error observed at different times is also qualitatively similar for IMEX and semi-implicit (see \Cref{fig:convergence_larsen}). This is likely because the adaptive time stepping for both methods drive the timestep sufficiently small that the IMEX method is locally in the basin of second-order convergence at any point in time and obtains better error than fixing a uniform timestep across the whole simulation.

\begin{table}[!htb]
\small
\centering
\begin{minipage}{\textwidth}
\centering
\begin{tabular}{c|ccc|ccc}
\hline
 & \multicolumn{3}{c|}{IMEX-NPRK2[42]b} & \multicolumn{3}{c}{H-LDIRK2(2,2,2)} \\
$\Delta t$ & $L_1$ & $L_2$ & $L_\infty$ & $L_1$ & $L_2$ & $L_\infty$ \\
\hline
3.20e-11 & 7.03e-02 & 9.90e-02 & 5.17e-01 & 7.18e-02 & 1.01e-01 & 5.22e-01 \\
1.60e-11 & 2.53e-02 (1.5) & 3.19e-02 (1.6) & 1.26e-01 (2.0) & 2.56e-02 (1.5) & 3.24e-02 (1.6) & 1.33e-01 (2.0) \\
8.00e-12 & 9.73e-03 (1.4) & 1.24e-02 (1.4) & 4.22e-02 (1.6) & 9.90e-03 (1.4) & 1.26e-02 (1.4) & 4.40e-02 (1.6) \\
4.00e-12 & 4.19e-03 (1.2) & 5.51e-03 (1.2) & 1.74e-02 (1.3) & 4.24e-03 (1.2) & 5.57e-03 (1.2) & 1.79e-02 (1.3) \\
2.00e-12 & 1.86e-03 (1.2) & 2.51e-03 (1.1) & 7.89e-03 (1.1) & 1.87e-03 (1.2) & 2.53e-03 (1.1) & 8.05e-03 (1.2) \\
1.00e-12 & 8.05e-04 (1.2) & 1.11e-03 (1.2) & 3.50e-03 (1.2) & 8.10e-04 (1.2) & 1.11e-03 (1.2) & 3.55e-03 (1.2) \\
\hline
\end{tabular}
\end{minipage}

\vspace{0.3cm}

\begin{minipage}{\textwidth}
\centering
\begin{tabular}{c|ccc|ccc}
\hline
 & \multicolumn{3}{c|}{SSP-LDIRK2(3,3,2)} & \multicolumn{3}{c}{SSP-LDIRK3(3,3,2)} \\
$\Delta t$ & $L_1$ & $L_2$ & $L_\infty$ & $L_1$ & $L_2$ & $L_\infty$ \\
\hline
3.20e-11 & 4.14e-02 & 5.70e-02 & 3.03e-01 & 2.93e-02 & 3.86e-02 & 2.05e-01 \\
1.60e-11 & 1.39e-02 (1.6) & 1.71e-02 (1.7) & 7.19e-02 (2.1) & 8.57e-03 (1.8) & 1.03e-02 (1.9) & 3.65e-02 (2.5) \\
8.00e-12 & 5.23e-03 (1.4) & 6.43e-03 (1.4) & 2.23e-02 (1.7) & 3.09e-03 (1.5) & 3.71e-03 (1.5) & 1.15e-02 (1.7) \\
4.00e-12 & 2.19e-03 (1.3) & 2.77e-03 (1.2) & 8.73e-03 (1.3) & 1.28e-03 (1.3) & 1.59e-03 (1.2) & 4.98e-03 (1.2) \\
2.00e-12 & 9.49e-04 (1.2) & 1.24e-03 (1.2) & 3.86e-03 (1.2) & 5.22e-04 (1.3) & 6.76e-04 (1.2) & 2.30e-03 (1.1) \\
1.00e-12 & 3.99e-04 (1.2) & 5.32e-04 (1.2) & 1.68e-03 (1.2) & 1.87e-04 (1.5) & 2.49e-04 (1.4) & 9.40e-04 (1.3) \\
\hline
\end{tabular}
\end{minipage}
\caption{\color{black}Fixed $\Delta t$ convergence for the temperature with IMEX integrators. Convergence rates are given in parenthesis with relative norm errors.}
\label{tab:convergence_imex_temperature}
\end{table}

\begin{table}[!htb]
\small
\centering
\begin{minipage}{\textwidth}
\centering
\begin{tabular}{c|ccc|ccc}
\hline
 & \multicolumn{3}{c|}{IMEX-NPRK2[42]b} & \multicolumn{3}{c}{H-LDIRK2(2,2,2)} \\
$\Delta t$ & $L_1$ & $L_2$ & $L_\infty$ & $L_1$ & $L_2$ & $L_\infty$ \\
\hline
3.20e-11 & 5.56e-02 & 1.59e-01 & 1.21e+00 & 5.58e-02 & 1.60e-01 & 1.22e+00 \\
1.60e-11 & 1.06e-02 (2.4) & 2.72e-02 (2.5) & 2.05e-01 (2.6) & 1.12e-02 (2.3) & 2.90e-02 (2.5) & 2.19e-01 (2.5) \\
8.00e-12 & 4.04e-03 (1.4) & 6.84e-03 (2.0) & 3.41e-02 (2.6) & 4.13e-03 (1.4) & 7.05e-03 (2.0) & 3.80e-02 (2.5) \\
4.00e-12 & 1.76e-03 (1.2) & 2.97e-03 (1.2) & 8.54e-03 (2.0) & 1.78e-03 (1.2) & 2.98e-03 (1.2) & 9.63e-03 (2.0) \\
2.00e-12 & 8.13e-04 (1.1) & 1.40e-03 (1.1) & 4.14e-03 (1.0) & 8.19e-04 (1.1) & 1.40e-03 (1.1) & 4.10e-03 (1.2) \\
1.00e-12 & 3.72e-04 (1.1) & 6.56e-04 (1.1) & 1.99e-03 (1.1) & 3.73e-04 (1.1) & 6.55e-04 (1.1) & 1.97e-03 (1.1) \\
\hline
\end{tabular}
\end{minipage}

\vspace{0.3cm}

\begin{minipage}{\textwidth}
\centering
\begin{tabular}{c|ccc|ccc}
\hline
 & \multicolumn{3}{c|}{SSP-LDIRK2(3,3,2)} & \multicolumn{3}{c}{SSP-LDIRK3(3,3,2)} \\
$\Delta t$ & $L_1$ & $L_2$ & $L_\infty$ & $L_1$ & $L_2$ & $L_\infty$ \\
\hline
3.20e-11 & 3.16e-02 & 8.75e-02 & 6.52e-01 & 1.99e-02 & 5.51e-02 & 4.14e-01 \\
1.60e-11 & 7.11e-03 (2.2) & 1.77e-02 (2.3) & 1.32e-01 (2.3) & 3.61e-03 (2.5) & 6.80e-03 (3.0) & 2.63e-02 (4.0) \\
8.00e-12 & 2.11e-03 (1.8) & 5.02e-03 (1.8) & 3.73e-02 (1.8) & 4.04e-03 (-0.2) & 1.04e-02 (-0.6) & 4.41e-02 (-0.7) \\
4.00e-12 & 8.22e-04 (1.4) & 1.92e-03 (1.4) & 1.43e-02 (1.4) & 1.96e-03 (1.0) & 5.12e-03 (1.0) & 2.14e-02 (1.0) \\
2.00e-12 & 3.64e-04 (1.2) & 8.58e-04 (1.2) & 6.45e-03 (1.2) & 9.39e-04 (1.1) & 2.47e-03 (1.1) & 1.05e-02 (1.0) \\
1.00e-12 & 1.59e-04 (1.2) & 3.80e-04 (1.2) & 2.88e-03 (1.2) & 4.38e-04 (1.1) & 1.17e-03 (1.1) & 5.20e-03 (1.0) \\
\hline
\end{tabular}
\end{minipage}
\caption{\color{black}Fixed $\Delta t$ convergence for $E_r$ with IMEX integrators. Convergence rates are given in parenthesis with relative norm errors.}
\label{tab:convergence_imex_rade}
\end{table}

\begin{table}[!htb]
\small
\centering
\begin{minipage}{\textwidth}
\centering
\begin{tabular}{c|ccc|ccc}
\hline
 & \multicolumn{3}{c|}{IMEX-NPRK2[42]b} & \multicolumn{3}{c}{H-LDIRK2(2,2,2)} \\
$\Delta t$ & $L_1$ & $L_2$ & $L_\infty$ & $L_1$ & $L_2$ & $L_\infty$ \\
\hline
3.20e-11 & 2.07e-04 & 3.84e-04 & 2.44e-03 & 2.23e-04 & 4.24e-04 & 1.95e-03 \\
1.60e-11 & 3.68e-05 (2.5) & 6.02e-05 (2.7) & 3.19e-04 (2.9) & 3.80e-05 (2.6) & 5.90e-05 (2.8) & 2.25e-04 (3.1) \\
8.00e-12 & 1.46e-05 (1.3) & 2.29e-05 (1.4) & 7.22e-05 (2.1) & 1.35e-05 (1.5) & 2.16e-05 (1.4) & 6.35e-05 (1.8) \\
4.00e-12 & 5.28e-06 (1.5) & 7.86e-06 (1.5) & 2.27e-05 (1.7) & 4.43e-06 (1.6) & 7.55e-06 (1.5) & 2.84e-05 (1.2) \\
2.00e-12 & 2.15e-06 (1.3) & 2.77e-06 (1.5) & 7.50e-06 (1.6) & 1.29e-06 (1.8) & 2.24e-06 (1.8) & 7.62e-06 (1.9) \\
1.00e-12 & 7.53e-07 (1.5) & 9.91e-07 (1.5) & 2.00e-06 (1.9) & 3.95e-07 (1.7) & 6.03e-07 (1.9) & 2.05e-06 (1.9) \\
\hline
\end{tabular}
\end{minipage}

\vspace{0.3cm}

\begin{minipage}{\textwidth}
\centering
\begin{tabular}{c|ccc|ccc}
\hline
 & \multicolumn{3}{c|}{SSP-LDIRK2(3,3,2)} & \multicolumn{3}{c}{SSP-LDIRK3(3,3,2)} \\
$\Delta t$ & $L_1$ & $L_2$ & $L_\infty$ & $L_1$ & $L_2$ & $L_\infty$ \\
\hline
3.20e-11 & 2.71e-04 & 4.34e-04 & 2.02e-03 & 1.43e-04 & 2.42e-04 & 1.28e-03 \\
1.60e-11 & 5.00e-05 (2.4) & 7.71e-05 (2.5) & 2.34e-04 (3.1) & 3.30e-05 (2.1) & 5.63e-05 (2.1) & 3.13e-04 (2.0) \\
8.00e-12 & 1.70e-05 (1.6) & 2.65e-05 (1.5) & 8.17e-05 (1.5) & 1.25e-05 (1.4) & 2.14e-05 (1.4) & 1.05e-04 (1.6) \\
4.00e-12 & 6.05e-06 (1.5) & 9.72e-06 (1.4) & 2.95e-05 (1.5) & 4.14e-06 (1.6) & 6.87e-06 (1.6) & 2.74e-05 (1.9) \\
2.00e-12 & 1.61e-06 (1.9) & 2.76e-06 (1.8) & 9.20e-06 (1.7) & 1.51e-06 (1.5) & 2.21e-06 (1.6) & 6.96e-06 (2.0) \\
1.00e-12 & 5.62e-07 (1.5) & 8.40e-07 (1.7) & 2.57e-06 (1.8) & 7.13e-07 (1.1) & 9.25e-07 (1.3) & 2.39e-06 (1.5) \\
\hline
\end{tabular}
\end{minipage}
\caption{\color{black}Fixed $\Delta t$ convergence for the temperature with semi-implicit integrators. Convergence rates are given in parenthesis with relative norm errors.}
\label{tab:convergence_semi_temperature}
\end{table}

\begin{table}[!htb]
\small
\centering
\begin{minipage}{\textwidth}
\centering
\begin{tabular}{c|ccc|ccc}
\hline
 & \multicolumn{3}{c|}{IMEX-NPRK2[42]b} & \multicolumn{3}{c}{H-LDIRK2(2,2,2)} \\
$\Delta t$ & $L_1$ & $L_2$ & $L_\infty$ & $L_1$ & $L_2$ & $L_\infty$ \\
\hline
3.20e-11 & 2.95e-04 & 4.87e-04 & 2.49e-03 & 2.83e-04 & 4.47e-04 & 2.16e-03 \\
1.60e-11 & 6.97e-05 (2.1) & 1.37e-04 (1.8) & 7.96e-04 (1.6) & 6.47e-05 (2.1) & 1.23e-04 (1.9) & 6.06e-04 (1.8) \\
8.00e-12 & 1.66e-05 (2.1) & 7.49e-05 (0.9) & 7.38e-04 (0.1) & 1.38e-05 (2.2) & 5.92e-05 (1.1) & 5.88e-04 (0.0) \\
4.00e-12 & 4.10e-06 (2.0) & 1.87e-05 (2.0) & 1.79e-04 (2.0) & 3.49e-06 (2.0) & 1.69e-05 (1.8) & 1.59e-04 (1.9) \\
2.00e-12 & 8.11e-07 (2.3) & 3.40e-06 (2.5) & 3.25e-05 (2.5) & 6.39e-07 (2.5) & 3.07e-06 (2.5) & 3.03e-05 (2.4) \\
1.00e-12 & 3.50e-07 (1.2) & 1.19e-06 (1.5) & 1.07e-05 (1.6) & 1.67e-07 (1.9) & 6.83e-07 (2.2) & 6.68e-06 (2.2) \\
\hline
\end{tabular}
\end{minipage}

\vspace{0.3cm}

\begin{minipage}{\textwidth}
\centering
\begin{tabular}{c|ccc|ccc}
\hline
 & \multicolumn{3}{c|}{SSP-LDIRK2(3,3,2)} & \multicolumn{3}{c}{SSP-LDIRK3(3,3,2)} \\
$\Delta t$ & $L_1$ & $L_2$ & $L_\infty$ & $L_1$ & $L_2$ & $L_\infty$ \\
\hline
3.20e-11 & 2.82e-04 & 4.22e-04 & 1.69e-03 & 2.86e-04 & 4.71e-04 & 2.78e-03 \\
1.60e-11 & 6.61e-05 (2.1) & 1.24e-04 (1.8) & 6.08e-04 (1.5) & 6.82e-05 (2.1) & 1.36e-04 (1.8) & 7.17e-04 (2.0) \\
8.00e-12 & 6.80e-06 (3.3) & 1.33e-05 (3.2) & 8.04e-05 (2.9) & 1.50e-05 (2.2) & 7.10e-05 (0.9) & 7.06e-04 (0.0) \\
4.00e-12 & 1.55e-06 (2.1) & 3.31e-06 (2.0) & 3.04e-05 (1.4) & 2.97e-06 (2.3) & 1.56e-05 (2.2) & 1.49e-04 (2.2) \\
2.00e-12 & 3.87e-07 (2.0) & 8.51e-07 (2.0) & 6.27e-06 (2.3) & 8.16e-07 (1.9) & 3.96e-06 (2.0) & 3.82e-05 (2.0) \\
1.00e-12 & 1.52e-07 (1.3) & 4.04e-07 (1.1) & 2.40e-06 (1.4) & 2.80e-07 (1.5) & 1.20e-06 (1.7) & 1.12e-05 (1.8) \\
\hline
\end{tabular}
\end{minipage}
\caption{\color{black}Fixed $\Delta t$ convergence for $E_r$ with semi-implicit integrators. Convergence rates are given in parenthesis with relative norm errors.}
\label{tab:convergence_semi_rade}
\end{table}
}

\section{Conclusions}\label{sec:conc}

In this paper we propose a moment-based adaptive time integration framework for thermal radiation transport. We extend the recently proposed IMEX integration framework for gray TRT \cite{imex-trt} to the multifrequency setting and to include a semi-implicit option, which is implicit in all components except opacity. We then derive embedded methods for four IMEX Runge--Kutta schemes that are asymptotic preserving for hyperbolic equations with stiff relaxation terms, a property our previous work has indicated is important for stable and robust IMEX integration for TRT and radiation hydrodynamics. Using the embedded methods, we develop error estimates in temperature and radiation energy that require trivial overhead in computational cost and storage, and enable the use of adaptive time stepping algorithms. We apply the new methods to the gray tophat and multifrequency Larsen problems, and demonstrate that the adaptive algorithm can naturally change the time step over 4--5 orders of magnitude in a single simulation, ranging from the streaming regime to the thick diffusion limit.

Currently, we only consider standalone TRT problems. Without coupling to other physics, these problems have a relatively simply evolution of dynamical timescale, starting very fast as radiation is streaming into the domain, and slowing down as the material heats up, eventually approaching an optically thick diffusive equation with much slower dynamical timescale. Perhaps not surprisingly, the adaptive procedure drives the time step smaller when the flux of radiation reaches a new material, including at time $t=0$ or as it reaches a new material as in Larsen, as well as when navigating sharp corners in the domain as in tophat. Future work will involve coupling the IMEX TRT formulation from \cite{imex-trt} and further developed here with our work on IMEX radiation hydrodynamics \cite{rad-hydro} for adaptive high energy density physics simulations, where we expect more dynamic behavior and interaction of materials and radiation.

\section*{Acknowledgments}
This work was supported by the Laboratory Directed Research and Development program of Los Alamos National Laboratory under project number 20220174ER, and DOE Office of Advanced Scientific Computing Research Applied Mathematics program through Contract No. 89233218CNA000001. LA-UR-25-20758. This manuscript has been authored by Lawrence Livermore National Security, LLC under Contract No. DE-AC52-07NA27344 with the US. Department of Energy. The United States Government retains, and the publisher, by accepting the article for publication, acknowledges that the United States Government retains a non-exclusive, paid-up, irrevocable, world-wide license to publish or reproduce the published form of this manuscript, or allow others to do so, for United States Government purposes.

\appendix
\section{Stability Regions}
\label{app:stability}

\begin{figure}[!htb]
    \centering

    \includegraphics[width=0.5\textwidth]{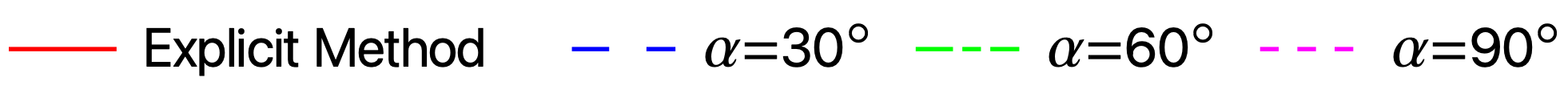}
    \\[0.75em]
    \begin{subfigure}{0.24\textwidth}
        \includegraphics[width=\textwidth]{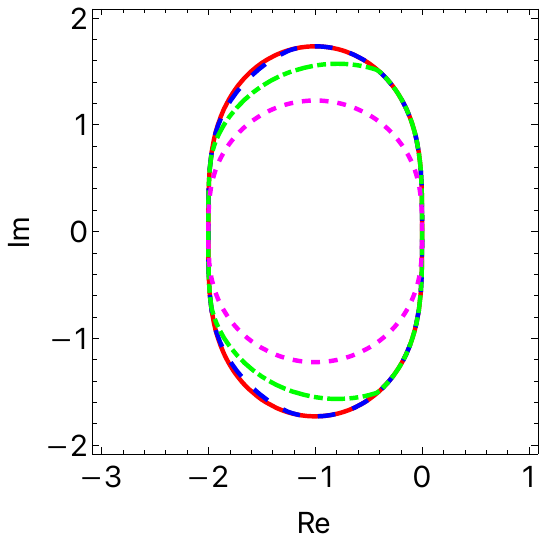}
        \caption{H-LDIRK2(2,2,2)}
    \end{subfigure}\hfill
    \begin{subfigure}{0.24\textwidth}
        \includegraphics[width=\textwidth]{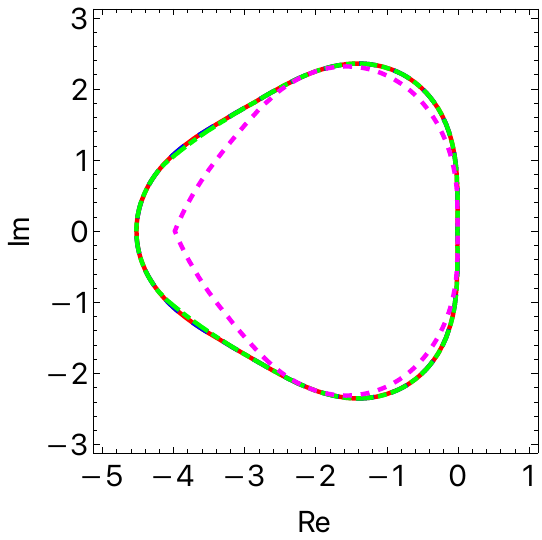}
        \caption{SSP-LDIRK2(3,3,2)}
    \end{subfigure}\hfill
    \begin{subfigure}{0.24\textwidth}
        \includegraphics[width=\textwidth]{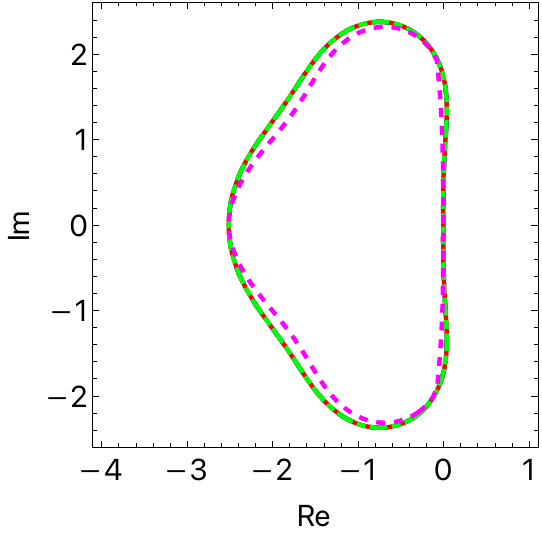}
        \caption{SSP-LDIRK3(3,3,2)}
    \end{subfigure}\hfill
    \begin{subfigure}{0.24\textwidth}
        \includegraphics[width=\textwidth]{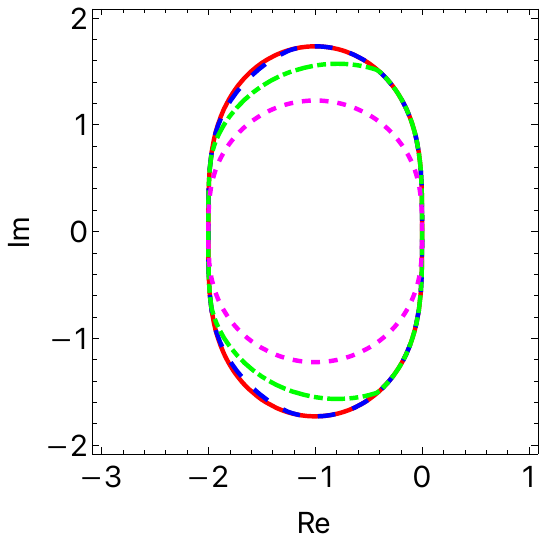}
        \caption{IMEX-NPRK2[42]b}
    \end{subfigure}
    \\[0.75em]
    \begin{subfigure}{0.24\textwidth}
        \includegraphics[width=\textwidth]{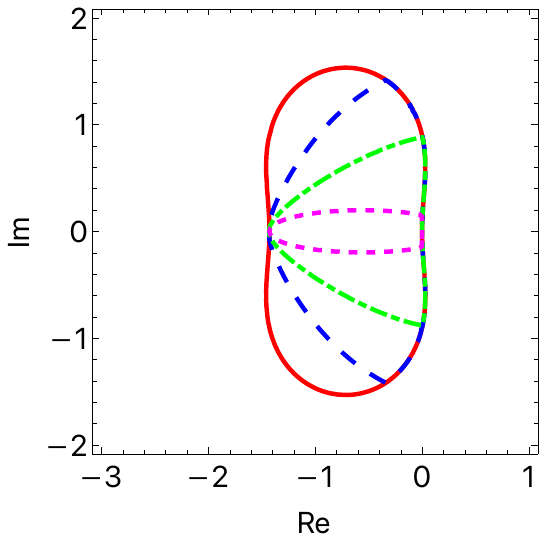}
        \caption{H-LDIRK2(2,2,2) Embedding}
    \end{subfigure}\hfill
    \begin{subfigure}{0.24\textwidth}
        \includegraphics[width=\textwidth]{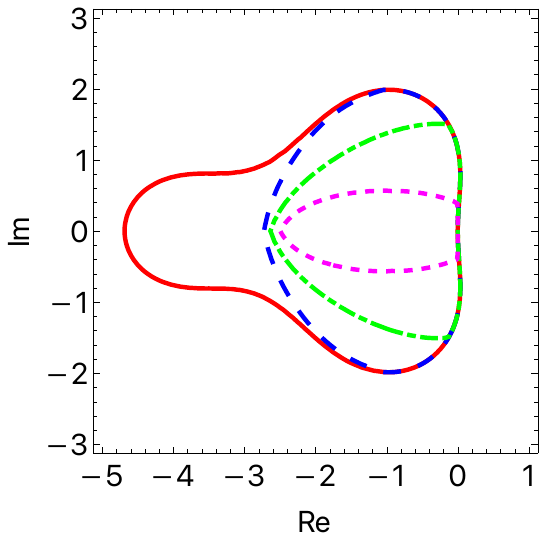}
        \caption{SSP-LDIRK2(3,3,2) Embedding}
    \end{subfigure}\hfill
    \begin{subfigure}{0.24\textwidth}
        \includegraphics[width=\textwidth]{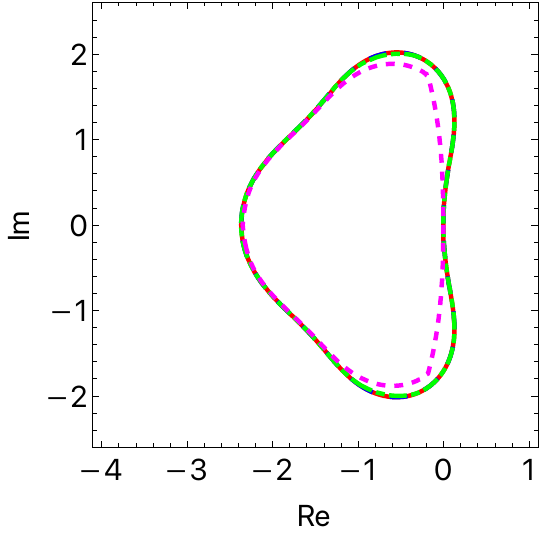}
        \caption{SSP-LDIRK3(3,3,2) Embedding}
    \end{subfigure}\hfill
    \begin{subfigure}{0.24\textwidth}
        \includegraphics[width=\textwidth]{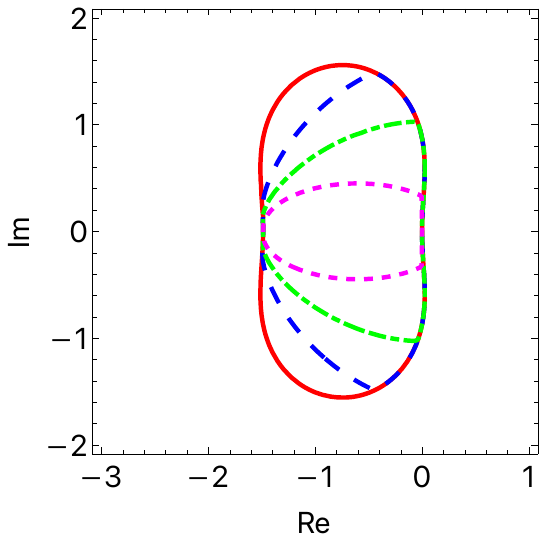}
        \caption{IMEX-NPRK2[42]b Embedding}
    \end{subfigure}
    
    \caption{The stability region $\mathcal{S}_{\infty, \alpha}^{\text{1D}}$ of the IMEX-RK schemes of \cref{sec:time:schemes} (top row) have the desirable property of being close in size to the stability region of the underlying explicit method.
    Corresponding stability regions for the embedded methods (bottom row) are slightly smaller.}
    \label{fig:stability}
\end{figure}

A common way to assess the linear stability of IMEX-RK schemes is by visualizing the region
\begin{equation} \label{eq:stability_region}
    \mathcal{S}_{\infty, \alpha}^{\text{1D}} = \left\{ \tilde{z} \in \mathbb{C} \middle| |R(\tilde{z}, z)| \leq 1, \forall z \in \mathbb{C}^{-} : |\arg(z) - \pi| \leq \alpha \right\}
\end{equation}
where $R(\tilde{z}, z)$ is the linear stability function associated with the test problem $y' = \tilde{\lambda} y + \lambda y$.
This shows the set of explicit scaled eigenvalues, $\tilde{z} = \Delta t \tilde{\lambda}$, such that for any implicit scaled eigenvalue, $z = \Delta t \lambda$, in a sector of angle $\alpha$, the method is stable.
We refer readers to \cite[Definition 4.1]{sandu2019class} for additional details.
\Cref{fig:stability} presents plots of \eqref{eq:stability_region} for the IMEX-RK methods of \cref{sec:time:schemes}.

\bibliographystyle{elsarticle-num} 
\bibliography{refs.bib}
\end{document}